\documentclass[11pt]{article}
\title{Isotopy invariants for closed braids and almost closed braids via loops in stratified spaces}
\author{Thomas Fiedler}
\usepackage{amssymb, amsfonts, epsfig}
\newtheorem{proposition}{Proposition}
\newtheorem{lemma}{Lemma}
\newtheorem{theorem}{Theorem}
\newtheorem{definition}{Definition}
\newtheorem{remark}{Remark}
\newtheorem{example}{Example}
\newtheorem{corollary}{Corollary}
\begin{document}
\maketitle
\begin{center}
{\em to S\'everine}
\end{center}
\begin{abstract}
Let $\phi : S^1\times D^2\to S^1$ be the natural projection. An oriented knot $K\hookrightarrow V = S^1\times D^2$ is called an almost closed braid 
if the restriction of $\phi$ to K has exactly two (non-degenerate) critical points (and K is a closed braid if the restriction of $\phi$ 
has no critical points at all).

We introduce new isotopy invariants for closed braids and almost closed braids in the solid torus V. These invariants refine finite type invariants.
They are still calculable with polynomial complexity with respect to the number of crossings of K.
Let the solid torus V be standardly embedded in the 3-sphere and let A be the axis of the complementary solid torus 
$S^3\setminus V$. We give examples which show that our invariants can detect non-invertibility of 2-component links 
$K\cup A\hookrightarrow S^3$. Notice that all quantum link invariants fail to do so and that it is not known wether there are finite type invariants 
which can detect non-invertibility of 2-component links.\footnote{2000 {\em Mathematics Subject Classification\/}: 57M25.{\em Key words and phrases\/}: Closed braids, Vassiliev invariants, character invariants, non-invertibility of links}
\end{abstract}
\tableofcontents
\section{Introduction}
This paper is a new and much shortened version of the preprint \cite{F03}.As well known, the isotopy problem for closed braids in the solid 
torus reduces to the conjugacy problem in braid groups (see \cite{M78}). The latter problem is solved, but in general only with exponential 
complexity with respect to the braid lenght (see \cite{B-B04} and references therein).
It is therefore interesting to construct invariants which distinguish conjugacy classes of braids and which are calculable 
with polynomial complexity. Finite type invariants for knots in the solid torus are an example of such invariants 
(see \cite{V90}, \cite{V01}, \cite{BN}, \cite{G96}, \cite{F01} and references therein).

In this paper we construct another class of calculable invariants.We give now a brief outline of our approach:  let $\hat \beta\subset V$ be a closed braid (i.e. the
restriction of $\phi$ to $\hat \beta$ has no critical points) and such that $\hat \beta$ is a knot. We fix a projection $pr : V \to S^1\times I$. Let $M(\hat \beta)$ be the  infinite dimensional) space of all closed braids which are isotopic to $\hat \beta$ in V. $M(\hat \beta)$ 
has a natural stratification. The strata $\sum^{(1)}$ of codimension 1 are just the braid diagrams which have in the projection $pr$ either 
exactly one ordinary triple point or exactly one ordinary autotangency.We call the corresponding strata $\sum^{(1)}(tri)$ and respectively 
$\sum^{(1)}(tan)$.
First, we associate to a closed braid in a canonical way a loop in $M(\hat \beta)$, called the {\em canonical loop   \/} and denoted
$rot(\hat \beta)$.Next, we associate to the canonical loop an oriented singular link in a thickened torus. This link is called the {\em trace graph
 \/}
and it is denoted by $TL(\hat \beta)$. All singularities of $TL(\hat \beta)$ are ordinary triple points. These triple points correspond exactly to the intersections of $rot(\hat \beta)$ with $\sum^{(1]}(tri)$.
There is a natural coorientation on $\sum^{(1)}(tri)$ and, hence, each triple point in $TL(\hat \beta)$ has a sign.
To each triple point corresponds a knot diagram which has just ordinary crossings and exactly one triple crossing.
We use the position of the ordinary crossings with respect to the triple crossing in the {\em Gauss diagram \/}
in order to construct {\em weights \/} for the triple points. We associate to the loop $rot(\hat \beta)$ then {\em weighted intersection numbers \/} with $\sum^{(1)}(tri)$. These intersection numbers turn out to be finite type invariants. We call them {\em one-cocycle invariants \/}.As all finite type invariants the one-cocycle invariants have a natural degree. However, they are rather special finite type invariants, as shows the following proposition.
\begin{proposition}
Let $\hat \beta$ be a closed n-braid which is a knot and let c be the word lenght of $\beta \in B_n$ (with respect to the standard generators of $B_n$).
Then all one-cocycle invariants of degree d vanish for $\hat \beta$  if    $d \geq  c + n^2 - n - 1$.
\end{proposition} 
One has to compare this proposition with the following well known fact:  the trefoil has non-trivial finite type invariants 
of arbitrary high degree.
Consequently, the one-cocycle invariants define a natural filtered subspace in the filtered space of all finite type invariants for those knots in the solid torus, which are closed braids.
There is a very simple procedure, coming from singularity theory, in order to construct all one-cocycle invariants for closed braids, without solving a big system of equations.
They verify automatically the marked 4T-relations (compare \cite{G96} and also \cite{F01}).

However , the one-cocycle invariants can be refined considerably.
When we deform $\hat \beta$ in V by a generic isotopy then $rot(\hat \beta)$ in $M(\hat \beta)$ deforms by a generic homotopy.
The following lemma is our key observation.
\begin{lemma}
The loop $rot(\hat \beta)$ is never tangential to $\sum^{(1)}(tan)$.
\end{lemma}

It follows from this lemma that the connected components of the natural resolution of $TL(\hat \beta)$ (i.e. the abstract union of circles 
where the branches in the triple points are separated) are isotopy invariants of $\hat \beta \hookrightarrow V$.
We apply now our theory of one-cocycle invariants but only to those triple points in $TL(\hat \beta)$ where three {\em given \/}
components of $TL(\hat \beta)$ intersect. The resulting invariants are called {\em character invariants \/}. They are no longer
of finite type but they are still calculable with polynomial complexity with respect to the braid lenght.

The set of finite type invariants (in particular , the set of one-cocycle invariants) can be seen as a trivial local system on 
$M(\hat \beta)$. In contrast to this, the set of character invariants is in general a non-trivial local system on $M(\hat \beta)$.

Alexander Stoimenow has written a computer program in order to calculate character invariants. It turns out that already 
character invariants of linear complexity can sometimes detect non-invertibility of closed braids  (i.e. the closed braid together with the axis of the complementary 
solid torus in $S^3$ is a non-invertible link in $S^3$). This implies in particular that character invariants are in general not of finite 
type.

We observe that the most simple character invariants are still well defined for almost closed braids.

.

The basic notions of our one parameter approach to knot theory, namely the space of non-singular knots, its discriminant,
the stratification of the discriminant, the coorientation of strata of low codimension, the canonical loop, the trace graph, the equivalence relation for trace graphs, are worked out in all 
details in our joint work with Vitaliy Kurlin \cite{F-K06}. Therefore we concentrate in this paper only on the construction of the new invariants.

\vspace{0,5cm}

{\em Acknowledgements \/}

I am grateful to Stepan Orevkov and Vitaliy Kurlin for many interesting discussions. I am especially grateful to Alexander Stoimenow for writing 
his computer program and for calculating interesting examples.
\section{Basic notions of one parameter knot theory}
In this section we recall briefly the basic notions of our theory. All details with complete proofs can be found in \cite{F-K06}
(even in a much more general setting).
\subsection{ The space of closed braids and its discriminant }
We work in the smooth category and all orientable manifolds are actually oriented.
We fix once for all a coordinate system in $\mathbb{R}^3$ : $(\phi , \rho , z )$.
Here, $(\phi ,\rho ) \in S^1\times \mathbb{R}^+$ are polar coordinates of the plane  $\mathbb{R}^2 = \{ z = 0 \}$.
A closed n-braid $\hat \beta$ is a knot in the solid torus $V = \mathbb{R}^3 \setminus  z-axes$, such that 
$\phi : \hat \beta \to  S^1$ is non-singular and $[\hat \beta] = n \in H_1(V)$. Let $M(\hat \beta)$ be the infinit dimensional space of all closed braids (with respect to $\phi$) which are isotopic to $\hat \beta$
in V. Let $M_n$ be the union of all spaces $M(\hat \beta)$.
A well known theorem of Artin (see e.g. \cite{M78}) says that two closed braids in the solid torus are isotopic as links in the solid torus 
if and only if they are isotopic as closed braids. Therefore it is enough to consider only isotopies through closed braids.
Let $pr : \mathbb{R}^3 \setminus z-axes \to \mathbb{R}^2 \setminus 0$ be the canonical projection 
$(\phi ,\rho , z) \to (\phi ,\rho )$.Each closed braid is then represented by a knot diagram with respect to $pr$. A generic closed braid $\hat \beta$ has only ordinary double points 
as singularities of $pr(\hat \beta)$. Let $\sum$ be the discriminant in $M(\hat \beta)$ which consists of all non-generic 
diagrams of closed braids isotopic to $\hat \beta$. 

The discriminant $\sum$ has a natural stratification:  $\sum = \sum^{(1)} \cup \sum^{(2)} \cup ...$,
where $\sum^{(i)}$ are the union of all strata of codimension i in $M(\hat \beta)$.

\begin{theorem}
(Reidemeisters theorem for closed braids)

$\sum^{(1)} = \sum(tri) \cup \sum(tan)$,

where $\sum(tri)$ is the union of all strata which correspond to diagrams with exactly one ordinary triple point 
(besides ordinary double points)  and $\sum(tan)$ is the union of all strata which correspond to diagrams with exactly one 
ordinary autotangency.
\end{theorem}

In the sequel we need also the description of $\sum^{(2)}$.

\begin{theorem}
$\sum^{(2)} = \sum_1 \cup \sum_2 \cup \sum_3 \cup \sum_4$, 

where $\sum_1$ is the union of all strata which correspond to diagrams with exactly one  ordinary quadruple point,
$\sum_2$ is the union of all strata which correspond to diagrams with exactly one ordinary autotangency 
through which passes another branch transversally,
$\sum_3$ corresponds to the union of all strata of diagrams with an autotangency in an ordinary flex,
$\sum_4$ is the union of all transverse intersections of strata from $\sum^{(1)}$.
\end{theorem}
\subsection{The canonical loop}

We identify $\mathbb{R}^3 \setminus z-axes$ with the standard solid torus $V = S^1 \times D^2 \hookrightarrow \mathbb{R}^3 \setminus z-axes$.
We identify the core of V with the unit circle in $\mathbb{R}^2$.

Let $rot(V)$ denote the $S^1$-parameter family of diffeomorphismes of V which is defined in the following way:
we rotate the solid torus monotoneously and with constant speed around its core by the angle t , $t \in [0 ,2\pi]$, i.e.
all discs $( \phi = const ) \times D^2$ stay invariant and are rotated simoultaneously around their centre.

 Let $\hat \beta$ be a closed braid.

\begin{definition}
The {\em canonical loop \/}  $rot(\hat \beta) \in M(\hat \beta)$ is the oriented loop induced by $rot(V)$.
\end{definition}

Notice that the whole loop $rot(\hat \beta)$ is completely determined by an arbitrary point in it.

The following lemma is an immediat corollary of the definition of the canonical loop .
\begin{lemma}
Let $\hat \beta_s , s\in [0 ,1 ]$ , be an isotopy of closed braids in the solid torus. Then $rot(\hat \beta_s ) , s \in [0 ,1 ]$,
is a homotopy of loops in $M(\hat \beta )$.
\end{lemma}
Evidently, the canonical loop can be defined for an arbitrary link in V in exactly the same way. However, in the case of closed braids
we can give an alternative combinatorial definition, which makes concret calculations much easier.

Let $\Delta \in B_n$ be Garside's element, i.e. $\Delta^2$ is a generator of the centre of $B_n$ (see \cite{B74}).
Geometrically, $\Delta^2$ is the full twist of the n strings.
\begin{definition}
Let $\gamma \in B_n$ be a braid with closure isotopic to $\hat \beta$. Then the {\em combinatorial canonical loop \/}
$rot(\gamma)$ is defined by the following sequence of braids:

$\gamma \to \Delta\Delta^{-1}\gamma \to \Delta^{-1}\gamma\Delta \to \dots \to \Delta^{-1}\Delta\gamma' \to \gamma' 
\to \Delta\Delta^{-1}\gamma' \to \Delta^{-1}\gamma'\Delta \to \dots \to \Delta^{-1}\Delta\gamma \to 
\gamma$Here, the first arrow consists only of Reidemeister II moves, the second arrow is a cyclic permutation of the braid word 
(which corresponds to an isotopy of the braid diagram in the solid torus) and the following arrows consist of "pushing $\Delta$
monotoneously from the right to the left through the braid $\gamma$". We obtain a braid $\tilde \gamma$ and we start again.
\end{definition}We give below a precise definition in the case $n = 3$. The general case is a straightforward generalization which is left to the reader.
$\Delta = \sigma_1\sigma_2\sigma_1$ for $n = 3$. We have just to consider the following four cases:$\sigma_1\Delta = \sigma_1(\sigma_1\sigma_2\sigma_1) \to \sigma_1(\sigma_2\sigma_1\sigma_2) = \Delta\sigma_2$$\sigma_2\Delta = (\sigma_2\sigma_1\sigma_2)\sigma_1 \to (\sigma_1\sigma_2\sigma_1)\sigma_1 = \Delta\sigma_1$

$\sigma_1^{-1}\Delta = \sigma_1^{-1}(\sigma_1\sigma_2\sigma_1) \to \sigma_2\sigma_1 \to (\sigma_1\sigma_1^{-1})\sigma_2\sigma_1 \to \sigma_1(\sigma_2\sigma_1\sigma_2^{-1}) = \Delta\sigma_2^{-1}$

$\sigma_2^{-1}\Delta = \sigma_2^{-1}(\sigma_1\sigma_2\sigma_1) \to (\sigma_1\sigma_2\sigma_1^{-1})\sigma_1 \to \sigma_1\sigma_2 \to \sigma_1\sigma_2\sigma_1\sigma_1^{-1} = \Delta\sigma_1^{-1}$.

Notice, that the sequence is canonical in the case of a generator and almost canonical in the case of an inverse generator. Indeed,
we could replace the above sequence   $\sigma_1^{-1}\Delta \to \Delta\sigma_2^{-1}$ by

$\sigma_1^{-1}(\sigma_1\sigma_2\sigma_1) \to \sigma_2\sigma_1 \to \sigma_2\sigma_1\sigma_2\sigma_2^{-1} \to (\sigma_1\sigma_2\sigma_1)\sigma_2^{-1}$.

But it turns out that the corresponding canonical loops in $M(\hat \beta)$ differ just by a homotopy which passes once transversally through 
a stratum of $\sum_2^{(2)}$.

Let c be the word lenght of $\gamma$. Then we use exactly $2c(n-2)$ braid relations (or Reidemeister III moves) in the combinatorial canonical loop.
This means that the corresponding loop in $M(\hat \beta)$ cuts $\sum^{(1)}(tri)$ transversally in exactly $2c(n-2)$ points.

One easily sees that the combinatorial canonical loop $rot(\gamma)$ from Definition 1 is homotopic in $M(\hat \beta)$ without touching $\sum^{(1)}(tan)$ 
(i.e. we never make
 in the one parameter family a Reidemeister II move forwards and just after that the same move backwards)
to the geometrical canonical loop $rot(\hat \beta)$ from Definition 2.
\subsection{The trace graph}
The trace graph $TL(\hat \beta)$ is our main combinatorial object. It is an oriented singular link in a thickened torus.
All its singularities are ordinary triple points.

Let $\hat \beta_t , t \in S^1$, be the (oriented) family of closed braids corresponding to the canonical loop $rot(\hat \beta)$.
We assume that the loop $rot(\hat \beta)$ is a generic loop. Let $\{ p_1^{(t)}, p_2^{(t}, \dots , p_m^{(t)} \}$ be the set of double points of
$pr(\hat \beta_t) \subset S_\phi^1 \times \mathbb{R}_\rho$.
The union of all these crossings for all $t \in S^1$ forms a link $TL(\hat \beta) \subset (S_\phi^1 \times \mathbb{R}_\rho^+) \times S_t^1$ (i.e. we forget the coordinate $z(p_i^{(t)})$).
$TL(\hat \beta)$ is non-singular besides ordinary triple points which correspond exactly to the triple points in the family $pr(\hat \beta_t)$.
A generic point of $TL(\hat \beta)$ corresponds just to an ordinary crossing $p_i^{(t)}$ of some closed braid $\hat \beta_t$.
Let $t: TL(\hat \beta) \to S_t^1$ be the natural projection. We orient the set of all generic points in $TL(\hat \beta)$
(which is a disjoint union of embedded arcs) in such a way that the local mapping degree of t at $p_i^{(t)}$ is $+1$ if and only if $p_i^{(t)}$
is a positive crossing (i.e. it corresponds to a generator of $B_n$ , or equivalently , its {\em writhe \/} $w(p_i^{(t)}) = +1$).

The arcs of generic points come together in the triple points and in points corresponding to an ordinary autotangency in some $pr(\hat \beta)$.
But one easily sees that the above defined orientations fit together to define an orientation on the natural resolution $TL\tilde (\hat \beta)$ of $TL(\hat \beta)$ 
(compare also \cite{Fi01}). $TL\tilde(\hat \beta)$ is a union of oriented circles, called {\em trace circles \/}.We can attach {\em stickers \/} $i \in \{ 1, 2, \dots ,n-1 \}$ to the edges of $TL(\hat \beta)$ in the following way:
each edge of $TL(\hat \beta)$ corresponds to a letter in a braid word. Indeed, each generic point in an edge corresponds to an ordinary crossing of a braid 
projection and , hence , to some $\sigma_i$ or some $\sigma_i^{-1}$ . We attache to this edge the number i . The information about the exponent $+1$
or $-1$ is contained in the orientation of the edge.

We identify $H_1(V)$ with $\mathbb{Z}$ by sending the core of V to the generator $+1$. If $\hat \beta$ is a knot then we can attache to each trace circle a 
{\em homological marking \/} $a \in H_1(V)$ in the following way: let p be a crossing corresponding to a generic point in the trace circle.
We smooth p with respect to the orientation of the closed braid. The result is an oriented 2-component link. The component of this link which contains 
the undercross  which goes to the overcross at p is called $p^+$. We associate now to p the homology class $a = [p^+] \in H_1(V)$ (compare also \cite{F93}).
One easily sees that $a \in \{ 1 ,2 , \dots , n-1\}$ and that the class a does not depend on the choice of the generic point in the trace circle.
Indeed, the two crossings involved in a Reidemeister II move have the same homological marking and a Reidemeister III move 
does not change the homological marking of any of the three involved crossings.

Evidently, crossings with different homological markings belong to different trace circles. Surprisingly, the inverse is also true.
The following lemma is essentially due to Stepan Orevkov.

\begin{lemma}
Let $TL(\hat \beta)$ be the trace graph of the closure of a braid $\beta \in B_n$, such that $\hat \beta$ is a knot.
Then $TL\tilde (\hat \beta)$ splitts into exactly n-1 trace circles. They have pairwise different homological markings.
\end{lemma}
Consequently, the trace circles are characterised by their homological markings. Notice, that the set of homological markings is independent of the word lenght of the braid.
\subsection{A higher order Reidemeister theorem for trace graphs of closed braids }

\begin{definition}
A {\em trihedron \/} is a 1-dimensional subcomplex of $TL(\hat \beta)$ which is contractible in the thickened torus and which has the form as shown in Figure 1.
\end{definition}
\begin{definition}
A {\em tetrahedron \/} is a 1-dimensional subcomplex of $TL(\hat \beta)$ which is contractible in the thickened torus and which has the form as shown in Figure 2.
\end{definition}
\begin{definition}A {\em trihedron move \/} is shown in Figure 3.
\end{definition}
\begin{definition}
A {\em tetrahedron move \/} is shown in Figure 4.
\end{definition}
\begin{figure}[htbp]
\centering \psfig{file=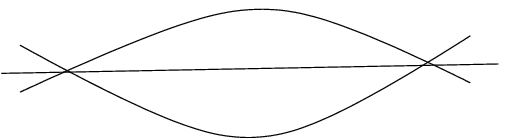}
\caption{}
\end{figure}
\begin{figure}[htbp]
\centering \psfig{file=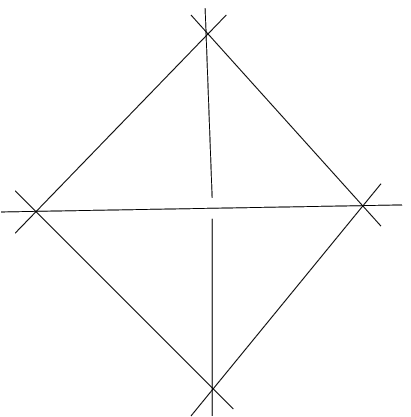}
\caption{}
\end{figure}
\begin{figure}[htbp]
\centering \psfig{file=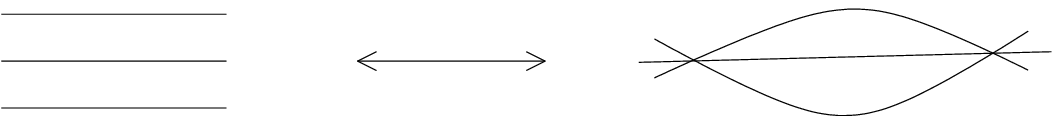}
\caption{}
\end{figure}
\begin{figure}[htbp]
\centering \psfig{file=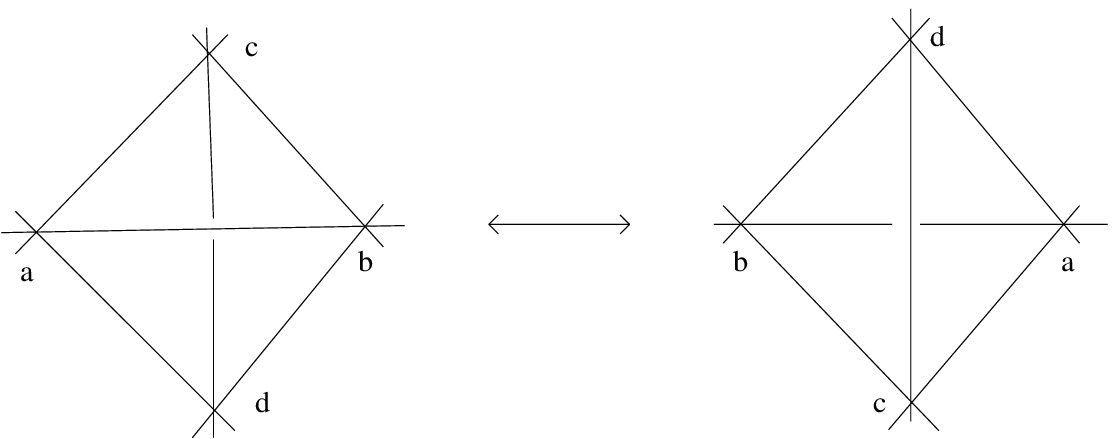}
\caption{}
\end{figure}
The rest of $TL(\hat \beta) \hookrightarrow S^1 \times S^1 \times \mathbb{R}^+$ is unchanged under the moves. The stickers on the edges change in the canonical way.

Notice, that a trihedron move corresponds to a generic homotopy of the canonical loop which passes once through an ordinary tangency 
with a stratum of $\sum^{(1)}(tri)$. A tetrahedron move corresponds to a generic homotopy of the canonical loop which passes transversally
once through a stratum of $\sum_1^{(2)}$, i.e. corresponding to an ordinary quadruple point.

\begin{definition}
The {\em equivalence relation for trace graphs $TL(\hat \beta)$ \/} is generated by the following three operations:

(1) isotopy in the thickened torus

(2) trihedron moves

(3) tetrahedron moves 
\end{definition}

The following important Reidemeister type theorem is a particular case of Theorem 1.10.  in \cite{F-K06}.

\begin{theorem}
Two closed braids (which are knots) are isotopic in the solid torus if and only if their trace graphs in the thickened torus are 
equivalent.
\end{theorem}

\begin{remark}
Notice, that not all representatives of an equivalence class of a trace graph correspond to the canonical loop (which is very rigid) of some closed braid. 
However, one easily sees that all representatives correspond to loops in the space $M(\hat \beta)$.
\end{remark}
\begin{remark}
A generic homotopy of a loop in $M(\hat \beta)$ could of course be tangent at some point to $\sum^{(1)}(tan)$. One easily sees that this would 
imply a Morse modification of the trace graph and, hence, change the components of its natural resolution. The corner stone of the theory developed in this paper is the fact,
that this does not happen for the (very rigid) homotopies of $rot(\hat \beta)$ which are induced by generic isotopies of $\hat \beta$ in V.
Consequently, the trace circles are isotopy invariants for closed braids (compare Lemma 3.4.  in \cite{F-K06}).
\end{remark}

\begin{definition}
A homotopy of loops $\gamma_s , s \in [0 , 1]$ , in $M(\hat \beta)$ is called a {\em tan-transverse homotopy \/} if no loop $\gamma_s$
is tangential to $\sum^{(1)}(tan)$. Let S be any loop in $M(\hat\beta)$. Its tan-transverse homotopy class is denoted by $[S]_{t-t}$.
\end{definition}

\begin{remark}
There can be loops $\gamma_s$ in a tan-transverse homotopy which are tangential to $\sum^{(1)}(tri)$. For the trace graphs associated to the loops 
$\gamma_s$ this corresponds to a trihedron move.
\end{remark}

In order to define our invariants we need more precise information about isotopies of trace graphs.

\begin{definition}
A {\em time section \/} in the thickened torus $(S^1_\phi \times \mathbb{R}^+_\rho) \times S^1_t$ is an annulus of the form 
$(S^1_\phi \times \mathbb{R}^+_\rho) \times \{ t = const \}$.
\end{definition}
The intersection of $TL(\hat \beta)$ with a generic time section corresponds to the crossings of the closed braid $\hat \beta_t$.
Using the orientation and the stickers on $TL(\hat \beta)$ we can read off a cyclic braid word for $\hat \beta$ in each generic time section.
The tangent points of $TL(\hat \beta)$  with  time sections correspond exactly to the Reidemeister II moves in the one parameter family of diagrams $\hat \beta_t , t \in S^1$.
A triple point in $TL(\hat \beta)$ slides over such a tangent point if and only if the canonical loop passes in a homotopy transversally through a stratum of $\sum_2^{(2)}$.
We illustrate this in Figure 5.

When the canonical loop passes transversally through a stratum of $\sum_3^{(2)}$ then the trace graph changes as shown in Figure 6.

Finally, when the canonical loop passes transversally through a stratum of $\sum_4^{(2)}$ then the t-values of triple points or tangencies with time sections 
are interchanged. We show an example in Figure 7.
\begin{figure}[htbp]
\centering \psfig{file=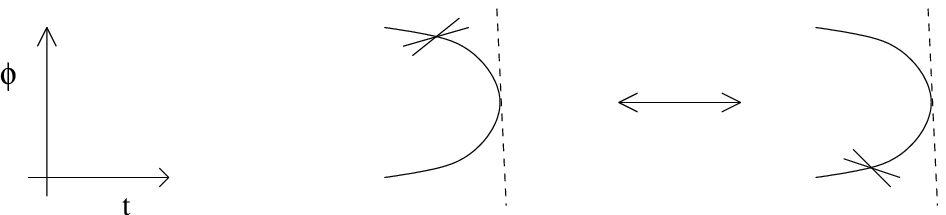}
\caption{}
\end{figure}
\begin{figure}[htbp]
\centering \psfig{file=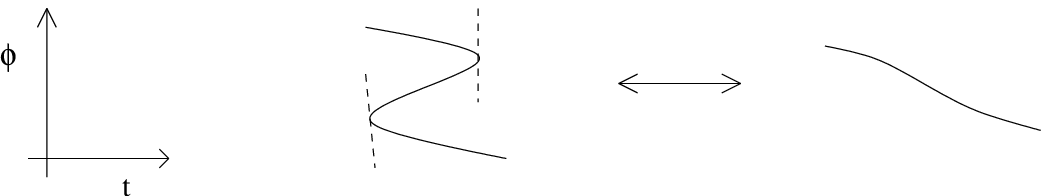}
\caption{}
\end{figure}
\begin{figure}[htbp]
\centering \psfig{file=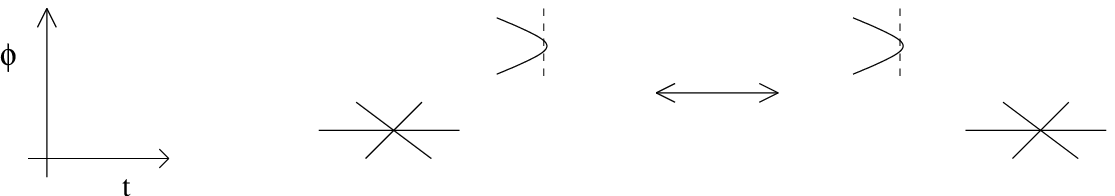}
\caption{}
\end{figure}
\section{One-cocycle invariants}

In this section we introduce our one-cocycle invariants. They are a special class of finite type invariants for knots in the solid torus.

\subsection{Gauss diagrams for closed braids with a triple crossing }
Let $f : S^1 \to \hat \beta$ be a generic orientation preserving diffeomorphisme. Let p be any crossing of $\hat \beta$.
We connect $f^{(-1)}(p) \in S^1$ by an oriented chord, which goes from the undercross to the overcross and we decorate it by the writhe $w(p)$.
Moreover, we attache to the chord the homological marking $[p^+]$. The result is called a {\em Gauss diagram \/} for
$\hat \beta$ (compare e.g. \cite{PV} and \cite{F01}).

One easily sees that $\hat \beta$ up to isotopy is determined by its Gauss diagram and the number $n = [\hat \beta] \in H_1(V)$.
(Notice, that the trace graph of a closed braid is always transverse to the $\phi$-sections.)

A {\em Gauss sum of degree d \/} is an expression assigned to a diagram of a closed braid which is of the following form:

$\sum$ function( writhes of the crossings)

where the sum is taken over all possible choices of d (unordered) different crossings in the knot diagram such that the chords without the writhes arising from these crossings 
build a given subdiagram with given homological markings. The marked subdiagrams (without the writhes) are called {\em configurations\/}.
If the function is the product of the writhes, then we will denote the sum shortly by the configuration itself.
We need to define Gauss diagrams for knots with an ordinary triple point too. The triple point corresponds to a triangle in the Gauss diagram of the knot.
Notice, that the preimage of a triple point has a natural ordering coming from the orientation of the $\mathbb{R}^+$-factor.
One easily sees that this order is completely determined by the arrows in the triangle.

We provide each stratum of $\sum^{(1)}(tri)$ with a coorientation which depends only on the non-oriented  underlying curves $pr(\hat \beta)$ in $S^1 \times \mathbb{R}^+$.
Consequently, for the definition of the coorientation we can replace the arrows in the Gauss diagram simply by chords.

\begin{definition}
The {\em coorientation \/} of the strata in $\sum^{(1)}(tri)$ is given in Figure 8. 
\end{definition}
Notice that the second line in Figure 8 does not occure for closed braids, but it does occure for almost closed braids.
The two coorientations are choosen in such a way that they fit together in $\sum^{(2)}_2$.

There are exactly two types of triple points without markings. We show them in Figure 9.

We attache now the homological markings to the three chords. Let $a,b \in \{1 , 2 , \dots , n-1\}$ be fixed. 
Then the markings are as shown in Figure 10. We encode the types of the marked triple points by $(a ,b)^-$ and $(a ,b)^+$.
The union of the corresponding strata of $\sum^{(1)}(tri)$ are encoded in the same way.
\begin{figure}[htbp]
\centering \psfig{file=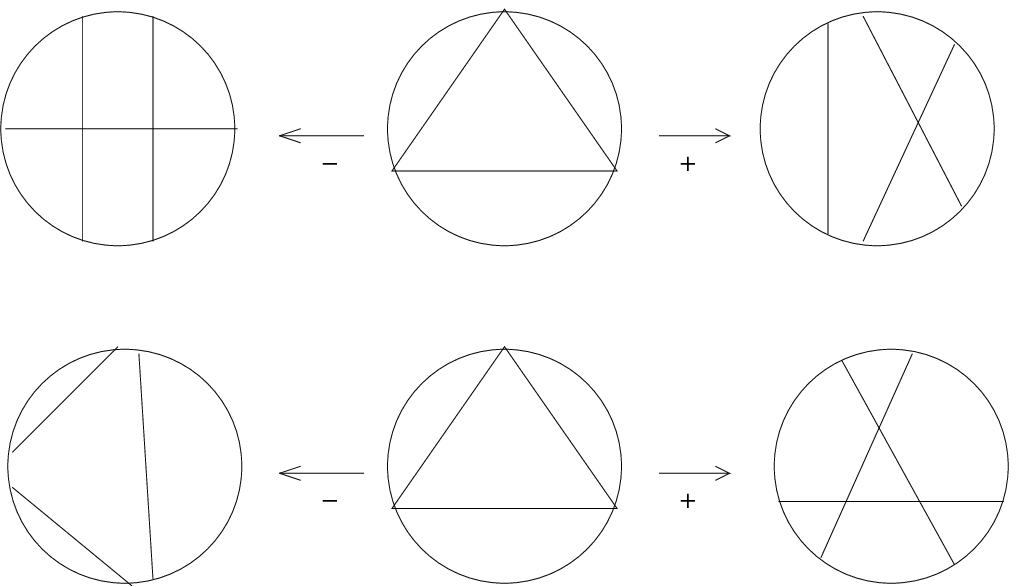}
\caption{}
\end{figure}
\begin{figure}[htbp]
\centering \psfig{file=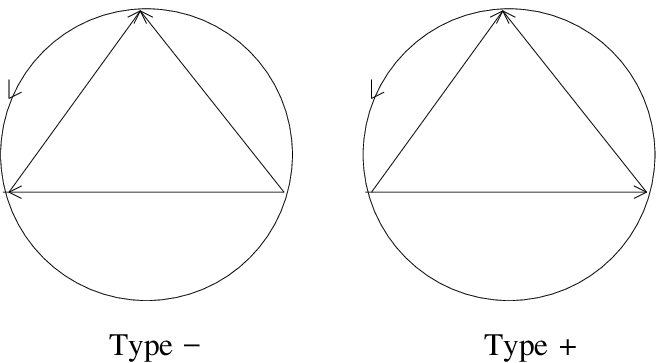}
\caption{}
\end{figure}
\begin{figure}[htbp]
\centering \psfig{file=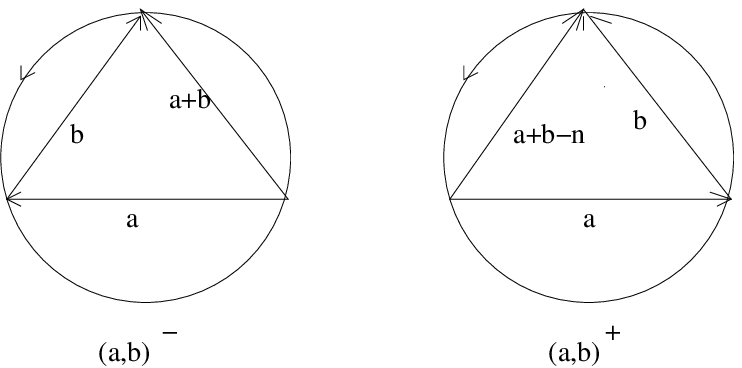}
\caption{}
\end{figure}
\subsection{One-cocycle invariants of degree one}

We will construct in a canonical way one-cocycles on the space $M_n$ (the space of all closed n-braids which are knots).
We obtain numerical invariants when we evaluate these cocycles on the homology class represented by the canonical loop.
 
Let $n, d \in \mathbb{N}^*$ be fixed. Let $(a,b)^\pm$ be a fixed
{\em type of marked triple point\/} as shown in Figure 10.
Here $a, b, a+b \in \mathbb{Z}/n\mathbb{Z}$.

\begin{definition}
A {\em configuration I of degree d\/} is an abstract Gauss diagram
without writhes which contains exactly $d-1$ arrows marked in
$\mathbb{Z}/n\mathbb{Z}$ besides the sub-diagram $(a,b)^\pm$.
\end{definition}

Let $\{ I_i \}$ be the finite set of all configurations of degree $d$
with respect to $(a,b)^\pm$. Let $\Gamma_{(a,b)^\pm} = \sum_{i}
{\epsilon_i I_i}$ be a linear combination with each $\epsilon_i \in
\{ 0, +1, -1 \}$. (The type of the triple point is {\em always\/}
fixed in any cochain $\Gamma$).
\begin{definition}
$\Gamma_{(a,b)^\pm}$ gives rise to a {\em 1-cochain of degree $d$\/}
by assigning to each oriented generic loop $S \subset M_n$
an integer $\Gamma_{(a,b)^\pm}(S)$ in the following way:
\begin{displaymath}
\Gamma_{(a,b)^\pm}(S) =
\sum_{s_i \in S \cap \Sigma^{(1)}(tri) \textrm{of type $(a,b)^\pm$}}
 {w(s_i)(\sum_{i}{\epsilon_i (\sum_{D_i}{\prod_{j}{w(p_j)}))}}}
\end{displaymath}
where $D_i$ is the set of unordered $(d-1)$-tuples $(p_1, \dots, p_{d-1})$
of arrows which enter in $I_i$ in the Gauss diagram of $s_i$ .
\end{definition}

\begin{lemma}
 If $\Gamma_{(a,b)^\pm}(S)$ is invariant under each generic deformation
of $S$ through a stratum of $\Sigma^{(2)}$, then $\Gamma_{(a,b)\pm}$
is a 1-cocycle.
\end{lemma}
{\em Proof:\/}  In this case, $\Gamma_{(a,b)^\pm}(S)$ is
invariant under homotopy of $S$. Indeed, tangent points of $S$ with $\Sigma^{(1)}(tri)$
correspond just to trihedron moves. The two triple points give the same contribution to 
$\Gamma_{(a, b)^\pm}(S)$ but with different signs. A tangency with $\Sigma^{(1)}(tan)$ does not change the contribution
of the triple points at all. This implies the invariance under
homology of $S$  because the
contributions to $\Gamma_{(a,b)\pm}(S)$ of different diagrams with
triple points are not related to each other. $\Box$

\begin{definition}
A {\em cohomology class\/} in $H^1(M_n; \mathbb{Z})$ is {\em
of degree $d$\/} if it can be represented by some 1-cocycle
$\Gamma_{(a,b)\pm}$ of degree at most $d$.
\end{definition}
\begin{remark}
The above definition induces a filtration on a part of $H^1(M_n;
\mathbb{Z})$.
\end{remark}

Let $M$ be the (disconnected) space of all embeddings $f: S^1
\hookrightarrow \mathbb{R}^3$. Vassiliev \cite{V90} has introduced a
filtration on a part of $H^1(M; \mathbb{Z})$ using the discriminant
$\Sigma_{sing}$. It is not difficult to see
that the space of all (unparametrized) differentiable maps of the
circle into the solid torus is contractible. Indeed, there is an
obvious canonical homotopy of each (perhaps singular) knot to a
multiple of the core of the solid torus.

The core of the solid torus is invariant under $rot_{S^1}(V)
\times rot_{D^2}(V)$. Thus, the above space is star-like. Therefore,
Alexander duality could be applied and Vassilievs approach could be
generalized for knots in the solid torus too.
It would be interesting to compare his filtration with our
filtration.

In the next sections, we will construct 1-cocycles $\Gamma_{(a,b)\pm}$
in an explicit way.

Let $\beta \in B_n$ be such that its closure $\hat \beta \hookrightarrow
V$ is a knot.

\begin{proposition}
 The space of finite type invariants of degree 1 is of dimension
$[n/2]$ (here $[.]$ is the integer part). It is generated by the
Gauss diagram invariants $W_a(\hat \beta) = \sum_{}{w(p)}$, where $a \in
\{ 1, 2, \dots, [n/2] \}$. The sum is over all crossings with fixed homological marking a.
\end{proposition}
{\em Proof:\/} It follows from Goryunov's \cite{G96} generalization of
finite type invariants for knots in the solid torus that the invariants
of degree 1 correspond just to marked chord diagrams with only one
chord. Obviously, all these invariants can be expressed as Gauss
diagram invariants:

$W_a(\hat \beta) = \sum_{}{w(p)}$, $a \in \{ 1, \dots, n-1 \}$

\vspace{0.5cm}

(see also \cite{F93}, and Section 2.2 in \cite{F01}.)

Let us define $V_a(\hat \beta) := W_a(\hat \beta) - W_{n-a}(\hat \beta)$
for all $a \in \{ 1, \dots, n-1 \}$. We observe that $V_a(\hat \beta)$
is invariant under switching crossings of $\hat \beta$.
Indeed, if the marking of the crossing $p$ was $[p] = a$, then the
switched crossing $p^{-1}$ has marking $[p^{-1}] = n-a$, but $w(p)=
-w(p^{-1})$.
But every braid $\beta \in B_n$ is homotopic to $\gamma = \prod_{i=1}
^{n-1}{\sigma_i}$. A direct calculation for $\gamma$ shows that
$V_a(\hat \gamma) \equiv 0$. It is easily shown by examples that
$W_a, a \in \{1, \dots, [n/2] \}$ (seen as invariants in $\mathbb{Q}$)
are linearly independent. $\Box$

\begin{lemma}
 Let $a, b \in \mathbb{Z}/n\mathbb{Z}$ be fixed.
Consider the union of all cooriented strata of $\Sigma^{(1)}$
which correspond to triple points of type either $(a, b)^-$ or $(a, b)^+$.
The closure in $M_n$ of each of these sets form 
integer cycles of codimension 1 in $M_n$.
\end{lemma}
\begin{remark}
Otherwise stated, $\Gamma_{(a,b)^+}$  and
$\Gamma_{(a,b)^-}$ both define integer 1-cocycles
of degree 1.
$\Gamma_{(a,b)^\pm}(S)$ is in this case by definition just the
algebraic intersection number of $S$ with the corresponding union
of strata of $\Sigma^{(1)}(tri)$.
\end{remark}

{\em Proof:\/} According to Section 2, we have to prove that the
cooriented strata fit together in $\Sigma^{(2)}_1$
and $\Sigma^{(2)}_2$. The first is evident , because at a stratum of $\Sigma^{(2)}_1$ just four strata of $\Sigma^{(1)}(tri)$ 
intersect pairwise transversally.
For the second, we have to distinguish 24 cases. Three of them are
illustrated in Figure 11.
The whole picture in a normal disc of $\Sigma^{(2)}_2$
is then obtained from Figure 12.

All other cases are obtained from these three by inverting the orientation
of the vertical branch, by taking the mirror image (i.e. switching
all crossings), and by choosing one of two possible closings of the
3-tangle (in order to obtain an oriented knot).
In all cases, one easily sees that the two adjacent triple points
are always of the same marked type and that the coorientations fit
together. $\Box$
\begin{figure}[htbp]
\centering \psfig{file=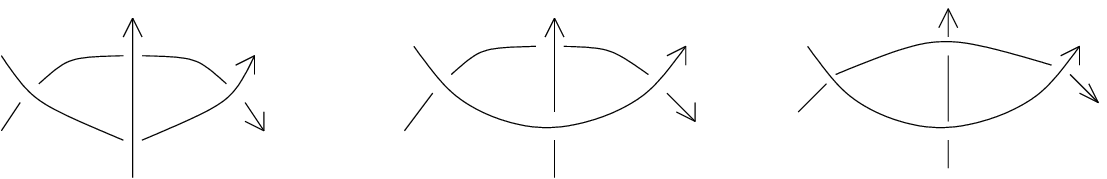}
\caption{}
\end{figure}
\begin{figure}[htbp]
\centering \psfig{file=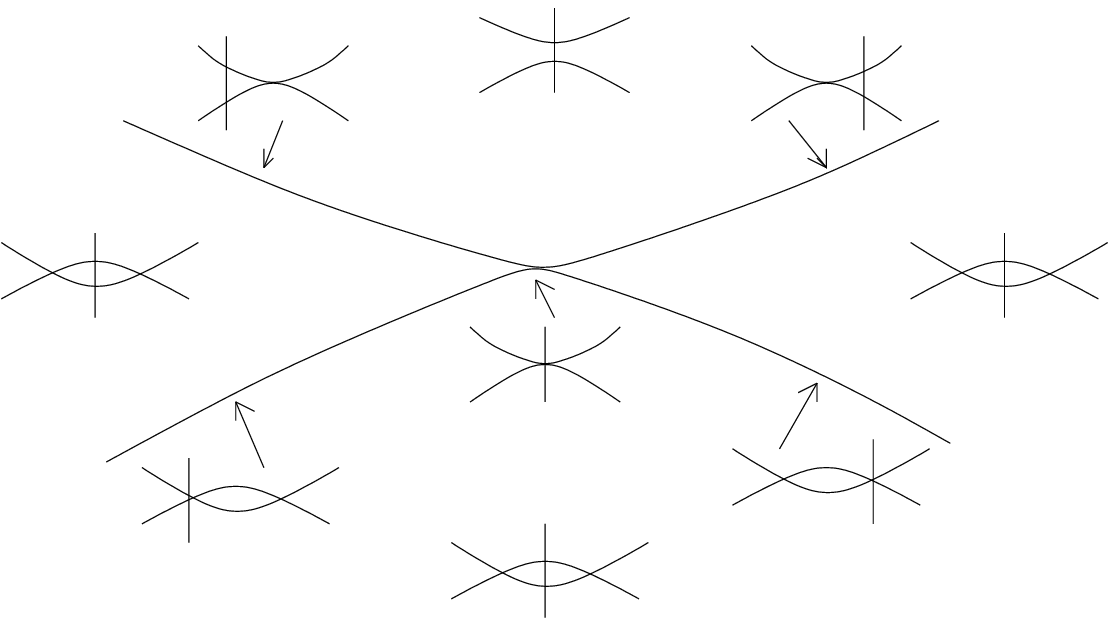}
\caption{}
\end{figure}
\begin{proposition}
$\Gamma_{(a,b)^+}$  defines a non-trivial 1-cohomology
class of degree 1 if and only if $a \not= b$ and $a + b \leq n-1$.

$\Gamma_{(a,b)^-}$  defines a non-trivial 1-cohomology
class of degree 1 if and only if $a \not= b$ and $a + b \geq n+1$.

The following identities hold:

(*) \hspace{4cm}    $\Gamma_{(a,b)^+} + \Gamma_{(b,a)^+} \equiv 0$

 (*) \hspace{4cm}   $\Gamma_{(a,b)^-} + \Gamma_{(b,a)^-} \equiv 0$

\end{proposition}

{\em Proof:\/} For closed braids, the markings are all in
$\{ 1, \dots, n-1 \}$. Therefore, if
$a+b>n-1$ in $\Gamma_{(a,b)^+}$ or $a+b<n+1$ in $\Gamma_{(a,b)^-}$,
then there is no such triple point at all and the 1-cocycle is trivial.
$\Gamma_{(a,a)^\pm} \equiv 0$ is a special case of the identities.

Examples show that all the remaining 1-cocycles are non-trivial
(we will see lots of examples later).

In order to prove the identities, we use the following Gauss
diagram sums (see also section 1.6 in \cite{F01}):
\begin{displaymath}
I^+_{(a,b)} = \sum_{} w(p)w(q), \qquad I^-_{(a,b)} = \sum_{} w(p)w(q)
\end{displaymath}
Here, the first sum is over all couples of crossings which form a subconfiguration as shown in Figure 13. The second sum is over all couples of 
crossings which form a subconfiguration as shown in Figure 14. (Here, a and b are the homological markings.)
These sums applied to diagrams of $\hat \beta$ are not invariants.
Let $S \subset M(\hat \beta)$ be a generic loop. Then $I^\pm_{(a,b)}$
is constant except when $S$ crosses $\Sigma^{(1)}(tri)$ in strata of type $(a,b)^\pm$ or
$(b,a)^\pm$.
At each such intersection in positive (resp., negative) direction,
$I^\pm_{(a,b)}$ changes exactly by $-1$ (resp., $+1$).
Indeed, the configurations of the three crossings which come together
in a triple point are shown in Figure 15.
In each of the four cases, exactly one pair $p, q$ of crossings contributes
to one of the sums $I$.

After drawing all possible triple points, it is easily seen that $p, q$
must verify: $w(p)w(q) = -1$ for the first two cases and
$w(p)w(q)=+1$ for the last two cases.
Notice that the type of the triple point is completely determined by
the sub-configuration shown in Figure 13 and 14.
Thus, the sums $I$ are constant by passing all types of triple
points except those shown in Figure 15. The generic loop $S$ intersects
$\Sigma$ only in strata that correspond to triple points or to
autotangencies. An autotangency adds to the Gauss diagram always one
of the sub-diagrams shown in Figure 16.
\begin{figure}[htbp]
\centering \psfig{file=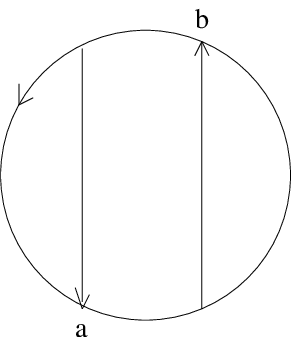}
\caption{}
\end{figure}
\begin{figure}[htbp]
\centering \psfig{file=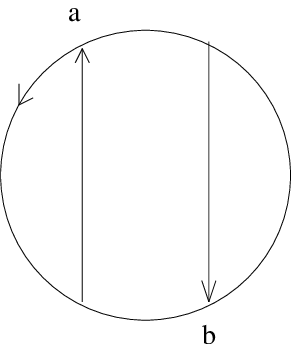}
\caption{}
\end{figure}
\begin{figure}[htbp]
\centering \psfig{file=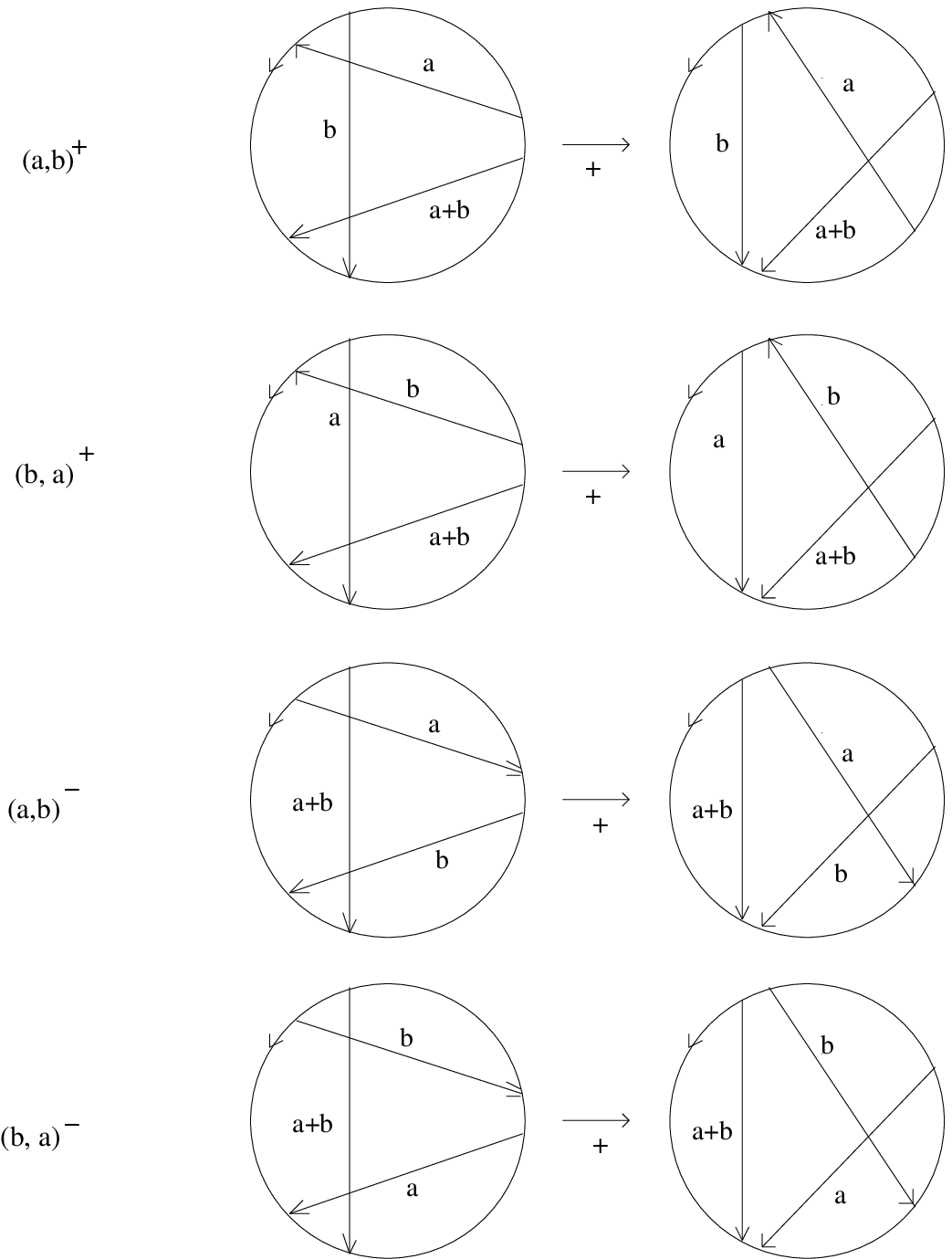}
\caption{}
\end{figure}
\begin{figure}[htbp]
\centering \psfig{file=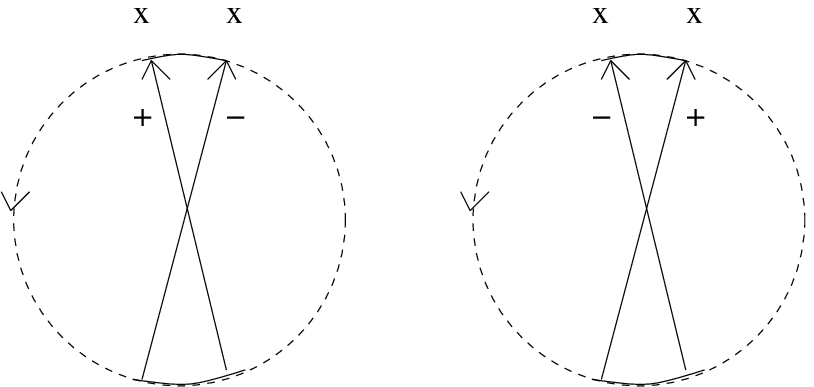}
\caption{}
\end{figure}
The two arrows do not enter together in the configurations shown in Figure 13 and 14. If one of them contributes to such a configuration,
then the other contributes to the same configuration but with an
opposite sign.

Therefore, for any $\hat \beta_1$, $\hat \beta_2 \in M(\hat \beta)
\setminus \Sigma$, the difference $I^\pm_{(a,b)}(\hat \beta_1) -
I^\pm_{(a,b)}(\hat \beta_2)$ is just the algebraic intersection
number of an oriented arc from $\hat \beta_1$ to $\hat \beta_2$
with the union of the cycles of codimension one
$(a,b)^+ \cup (b,a)^+$ (resp., $(a,b)^- \cup (b,a)^-$).
Hence, for each loop $S$, these numbers are 0. The identities (*) follow
and the proof of the proposition is complete.$\Box$

\begin{remark}
Obviously, for closed 2-braids, there are no 1-cocycles of any degree
$d$ at all (because there are never triple points). The previous proposition
implies that there are no non-trivial 1-cocycles of degree one for
closed 3-braids, but that in general there are $2((n-3)+(n-5)+(n-7)+
\dots)$ such cocycles.

We will see later that if we apply a 1-cocycle of degree one to the
canonical class $[rot(\hat \beta)]$ then we obtain a finite
type invariant of degree one for $\hat \beta$.
Proposition 3 says that there are exactly $[n/2]$ such invariants
which are independent. But the number of non-trivial 1-cocycles of
degree one is quadratic in $n$. Therefore, there are lots of
relations between them when they are restricted to the canonical
class. The following question seems to be interesting:
 Are the relations (*) from Proposition 3 the only relations
in general between the 1-cocycles of degree one?
\end{remark}
\begin{example}
Let $\hat \beta$ be the closure of the 4-braid $\beta = \sigma_1
\sigma_2^{-1}\sigma_3^{-1}$. We consider
\begin{displaymath}
\Gamma_{(1,2)^-} and \Gamma_{(2,1)^-}
\end{displaymath}
\begin{displaymath}
\Gamma_{(2,3)^+} and \Gamma_{(3,2)^+}
\end{displaymath}
A calculation by hand gives:
\begin{displaymath}
\Gamma_{(1,2)^-}(rot(\hat \beta)) =
\Gamma_{(2,3)^+}(rot(\hat \beta)) = -1
\end{displaymath}
and
\begin{displaymath}
\Gamma_{(2,1)^-}(rot(\hat \beta)) =
\Gamma_{(3,2)^+}(rot(\hat \beta)) = +1
\end{displaymath}
Therefore, all four 1-cocycles of degree 1 are non-trivial.
\end{example}
\subsection{One-cocycles of degree two}
Let $\beta \in B_n$ such that $\hat \beta \hookrightarrow V$ is a knot
and let $a \in \mathbb{Z}/n\mathbb{Z}$ be fixed. Let $I$ be a
configuration of degree 2, i.e. besides the marked triple point there
is exactly one marked arrow.

\begin{definition}
An {\em adjacent configuration of $I$\/} is any configuration which
is obtained by sliding exactly one of the end points of the arrow
over exactly one of the vertices of the triangle and which preserves the
markings. An example is shown in Figure 17.
A {\em chain of adjacent configurations\/} is a sequence of
configurations such that any two consecutive configurations are
adjacent.
\end{definition}
\begin{definition}
A configuration $I$ of degree 2 is called {\em braid impossible\/}
if it never occures as a subconfiguration of a Gauss diagram of a
closed $n$-braid which is a knot. Otherwise it is called {\em braid
possible\/}
\end{definition}
\begin{lemma}
 Any configuration $I$ which contains one of the sub-configurations
in Figure 18 is braid impossible.
\end{lemma}
{\em Proof:\/} This follows immediately from the proof of
Proposition 4.2 in \cite{Fi03}.$\Box$

\begin{remark}
There are other configurations $I$ which are braid impossible.
For example, it follows from Proposition 4.3 in \cite{Fi03} that $I$ is
braid impossible for $n=3$ if it contains one of the sub-configurations
drawn in Figure 19.
\end{remark}
\begin{figure}[htbp]
\centering \psfig{file=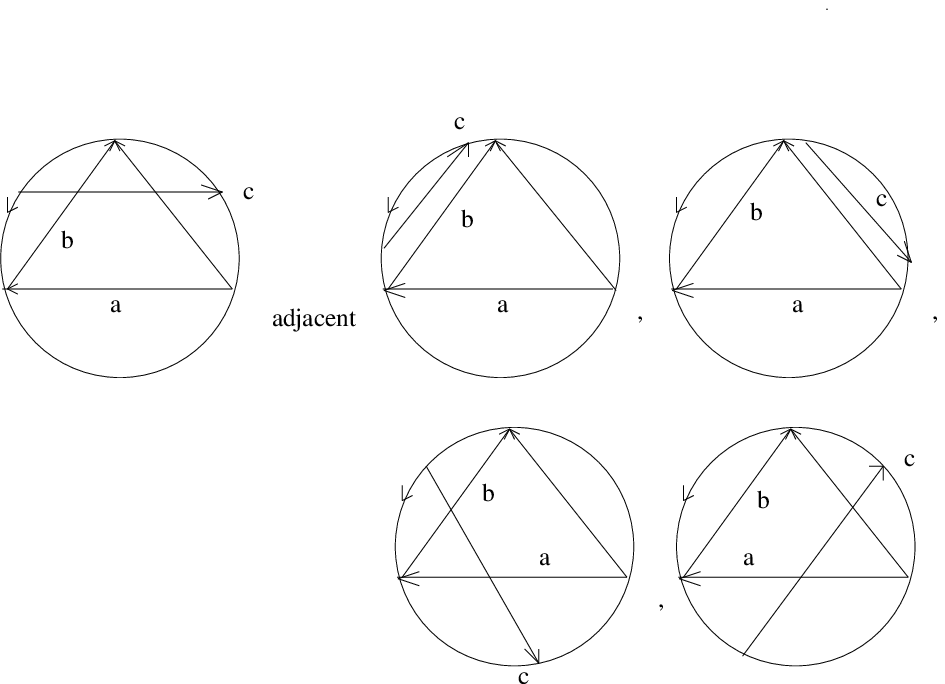}
\caption{}
\end{figure}
\begin{figure}[htbp]
\centering \psfig{file=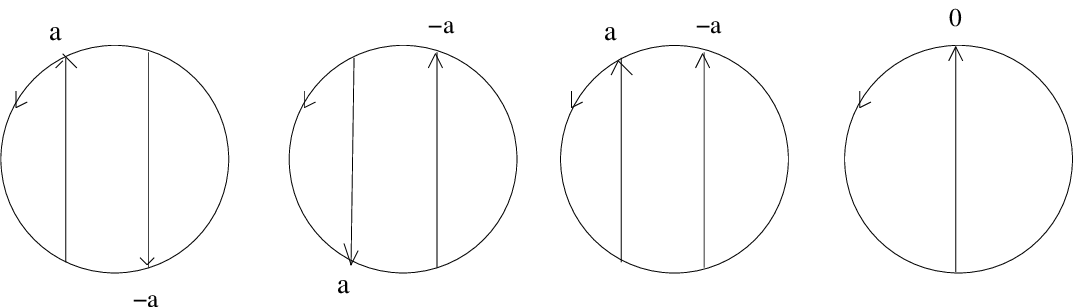}
\caption{}
\end{figure}
\begin{figure}[htbp]
\centering \psfig{file=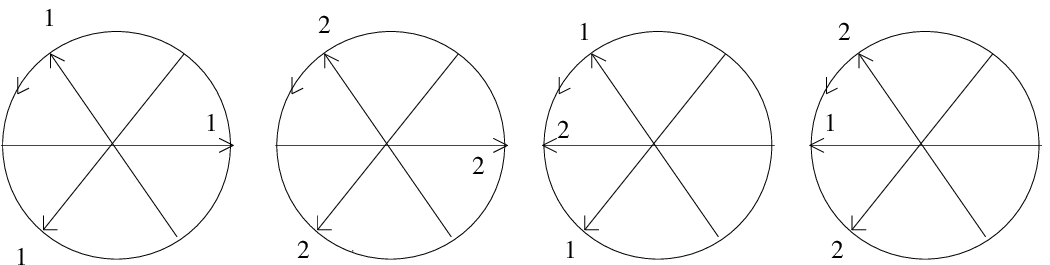}
\caption{}
\end{figure}
\begin{definition}
A configuration $I$ is called {\em rigid\/} if {\em all\/}
adjacent configurations are braid impossible. Otherwise it is called
{\em flexible\/}.
\end{definition} 

\begin{lemma}
 Let $a \in \{ 1, 2, \dots ,[(n-1)/2] \}$ be fixed.The configurations
shown in Figure 20 are all braid possible and rigid.
\end{lemma}
{\em Proof:\/} It is easy to show that the above configurations are braid
possible by just taking examples. Let us show that e.g. the first of
them is rigid. Indeed, all four adjacent configurations contain the
third sub-configuration of Figure 18. Thus, by Lemma 6, they are not
braid possible. The proof in the other cases is analogous.$\Box$
\begin{remark}
Let $I$ be a braid possible and rigid configuration. Let $I'$ be any
configuration obtained from $I$ by reversing some arrows and replacing for
each of these arrows its marking $x$ by $n-x$.
We say that $I'$ is {\em derived from\/} $I$.
If the triangle of $I'$ corresponds still to a Reidemeister III move, then $I'$ is also braid possible and rigid.
This is evident from the fact
that the set of closed braids is invariant under arbitrary changings of
crossings. Therefore, the above construction allows to get lots of
braid possible rigid configurations starting from the three cases in
Figure 20 (examples are shown in Figure 21).
\end{remark}
Let $\Gamma(S)$ be a 1-cochain of degree 2. Obviously, $\Gamma(S)$
is invariant under deformation of $S$ through all strata in
$\Sigma^{(2)}$ besides $\Sigma^{(2)}_1$ and $\Sigma^{(2)}_2$
(see the proof of Proposition 3).
Let us consider quadruple points. We have to guarantee that $\Gamma(S_m)=0$
if $S_m$ is the boundary of a small normal disc of a stratum in
$\Sigma^{(2)}_1$. Each triple point occurs exactly twice in $S_m$, say at
$s_1, s_2 \in S_m$., and with different signs $w(s_1)=-w(s_2)$.

The Gauss diagrams of the closed braids are the same at $s_1$ and
$s_2$ besides exactly three crossings $p_1, p_2, p_3$ which are now
in another position with respect to the triple point (and to each
other). We illustrate this in Figure 22.
Let $I$ be a configuration which contains the triple point and
only one of the three crossings, say $p_1$. One easily sees that the
triple point together with $p'_1$ is then an adjacent configuration.
For example, $p_1$ and $p'_1$ are in the same position with respect to
the crossing between the branches 2 and 3 (see Figure 22).
{\em Thus, if $n>3$, then for each configuration $\epsilon_iI_i$ in
$\Gamma$ (see Definition 12), $\Gamma$ has to contain all
adjacent configurations of $I_i$ and they have to enter all with the
same coefficient $\epsilon_i$.\/} We call this condition on $\Gamma$
the {\em quadruple-condition\/}. If $\Gamma$ verifies the
{\em quadruple-condition\/},
then $\Gamma(S_m)=0$. In particular, this implies that the triple
point together with $p_1$ or $p_2$ or $p_3$ is never a rigid
configuration! It remains to study a loop $S_m$ which is the boundary
of a normal disc for a stratum in $\Sigma^{(2)}_2$ (see Figure 12).
The two triple points in $S_m$ are of the same type and have different
signs as already explained in the proof of Lemma 5. But their
Gauss diagrams differ by exactly one crossing as shown in Figure 23.
But now $p$ and $p'$ are different crossings. The homological markings
$[p]=[p']$ coincide but $w(p)=-w(p')$.

Therefore, the configurations on the left-hand side and on the
right-hand side of Figure 23 should enter in $\Gamma$ with opposite
coefficients $\epsilon_i$. Notice that the marking $[p]$ of $p$
coincides with at least one of the markings of the triple point.
\begin{figure}[htbp]
\centering \psfig{file=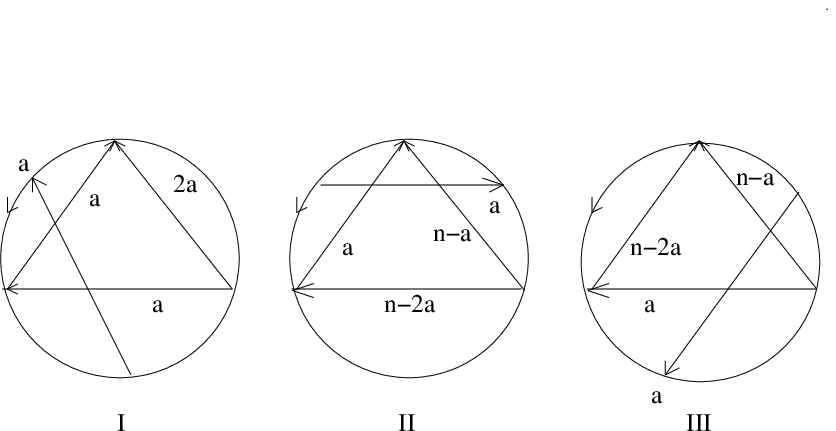}
\caption{}
\end{figure}
\begin{figure}[htbp]
\centering \psfig{file=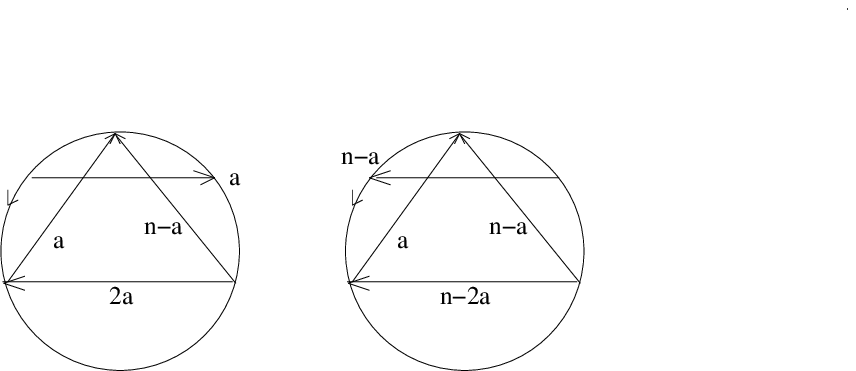}
\caption{}
\end{figure}
\begin{figure}[htbp]
\centering \psfig{file=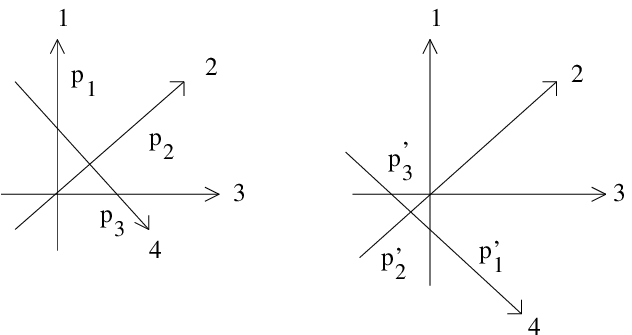}
\caption{}
\end{figure}
\begin{figure}[htbp]
\centering \psfig{file=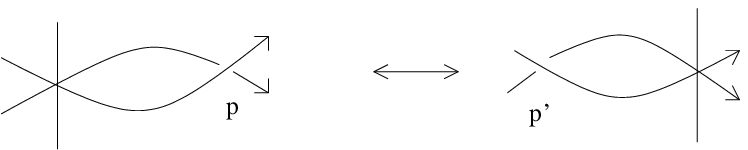}
\caption{}
\end{figure}
\begin{definition}
A configuration $I$ is called {\em t-invariant\/} if the
marking of the arrow is different from all three markings of the triple
point.
\end{definition}

Let $\Gamma' \subset \Gamma$ be the t-invariant
part of $\Gamma$. Then, evidently, $\Gamma'(S_m)=0$.

Let us now consider configurations which are not t-invariant.
If the two configurations in Figure 23 could be related by a chain of
adjacent configurations then our method would break down.
Indeed, in order to guarantee invariance under passing a quadruple point,
they should enter in $\Gamma$ with the same coefficient $\epsilon$.
But in order to guarantee invariance under passing an autotangency with a transverse branch they
should enter in $\Gamma$ with opposite coefficients. The following surprising
lemma implies that this does not occur for some types of triple
points.
\begin{lemma}
 Let the type of the triple point in Figure 23 be one of the 
types shown in Figure 20 (or one obtained from them as explained in
Remark 8). Then at least one of the two configurations in
Figure 23 is rigid.
\end{lemma}
{\em Proof:\/} We have to distinguish two cases for the markings.

{\em Case~1\/} All three markings of the triple point are different.
This is the case in $II$ and $III$ of Figure 20 if $3a \not= n$.

{\em Case~2\/} There are exactly two markings which are equal, hence
we are in $I$ of Figure 20 or in $II$ or $III$ with $3a=n$. (Remember
that markings of braid possible configurations are always non-zero).

The crossings $p$ and $p'$ form together a sub-configuration as shown
in Figure 18. The crossings $p$ and $p'$ interchange the place in the
triangle (corresponding to the triple point) and one easily sees that
both arrows $p$ and $p'$ always move in the same direction on the
circle. Therefore, in Case~1, we have the couples of configurations
as shown in Figure 24.
We see that the configurations on the left-hand side are always
rigid and those on the right-hand side are always flexible.
In Case~2, we have the couples of configurations shown in Figure 25.

\begin{figure}[htbp]
\centering \psfig{file=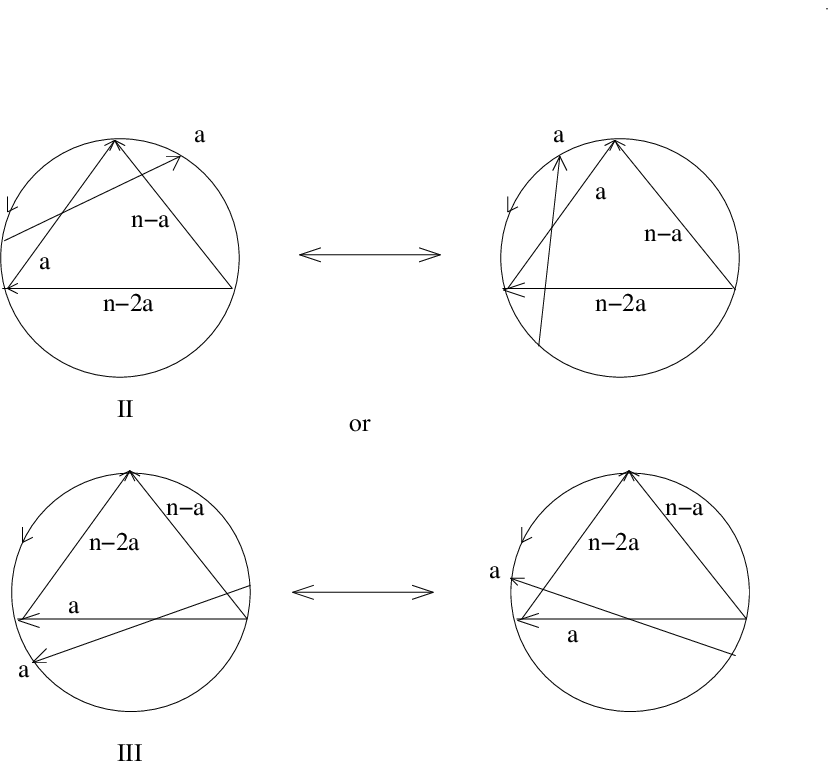}
\caption{}
\end{figure}
\begin{figure}[htbp]
\centering \psfig{file=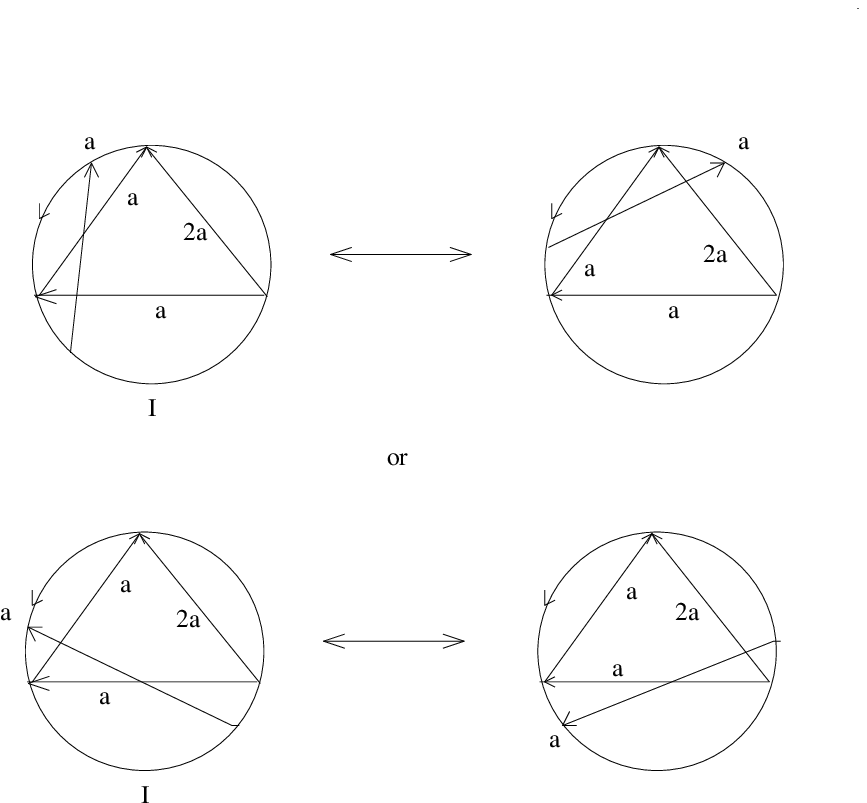}
\caption{}
\end{figure}
The cases $I$ and $III$ with $3a=n$ are equivalent to the cases shown
in Figure 25. We see that both configurations are always rigid.$\Box$

{\em Thus, for each configuration $\epsilon_iI_i$ in $\Gamma$ which
is shown on the right-hand side of Figure 24 or 25, $\Gamma$ has to
contain also the corresponding rigid configuration on the left-hand
side of Figure 24 or 25. This configuration has to enter in
$\Gamma$ with the coefficient $-\epsilon_i$.\/}
We call this condition on $\Gamma$ the {\em t-condition\/}.

\begin{theorem}
 Let $\Gamma$ be a 1-cochain of degree 2 that satisfies the quadruple-condition and the t-condition.
Then, $\Gamma$ is a 1-cocycle of degree 2.
\end{theorem}

{\em Proof:\/} We have proven that $\Gamma(S)$ is invariant under
all homotopies of $S$. The rest of the proof follows from Lemma  8.
$\Box$
\begin{example}
The case of closed 3-braids is very special because there do not appear
any quadruple points in isotopies. Therefore we do not need the
quadruple-condition. It follows easily from the proof of Lemma 8
that the 1-cochain in Figure 26 defines a 1-cocycle of degree 2.
(For closed 3-braids, the homological markings of the triple point are
determined by the arrows.) An easy calculation yields:
\begin{displaymath}
\Gamma(rot(\hat {\sigma_1\sigma_2^{-1}}))=+2
\end{displaymath}
Therefore, $\Gamma$ defines a non-trivial 1-cohomology class of degree 2.

It follows from Theorem 4 that the $\Gamma$ shown in Figure 27
defines a 1-cocycle of degree 2 for closed 4-braids.
A calculation yields:
\begin{displaymath}
\Gamma(rot(\hat {\sigma_1\sigma_2^{-1}\sigma_3^{-1}}))=-1
\end{displaymath}
Therefore, $\Gamma$ is a non-trivial class of degree 2.
\end{example}
\begin{figure}[htbp]
\centering \psfig{file=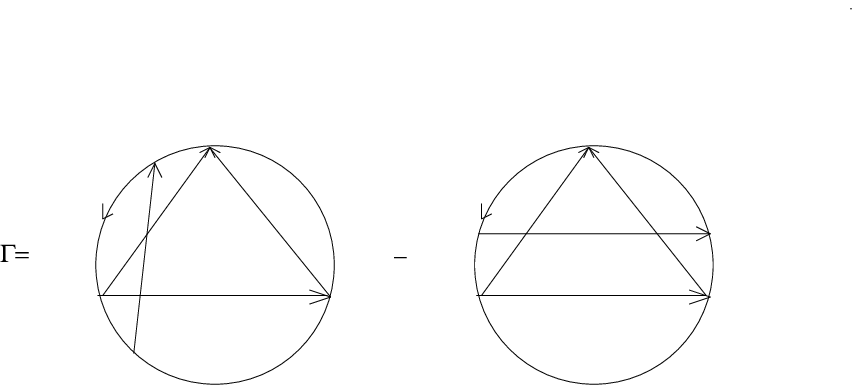}
\caption{}
\end{figure}
\begin{figure}[htbp]
\centering \psfig{file=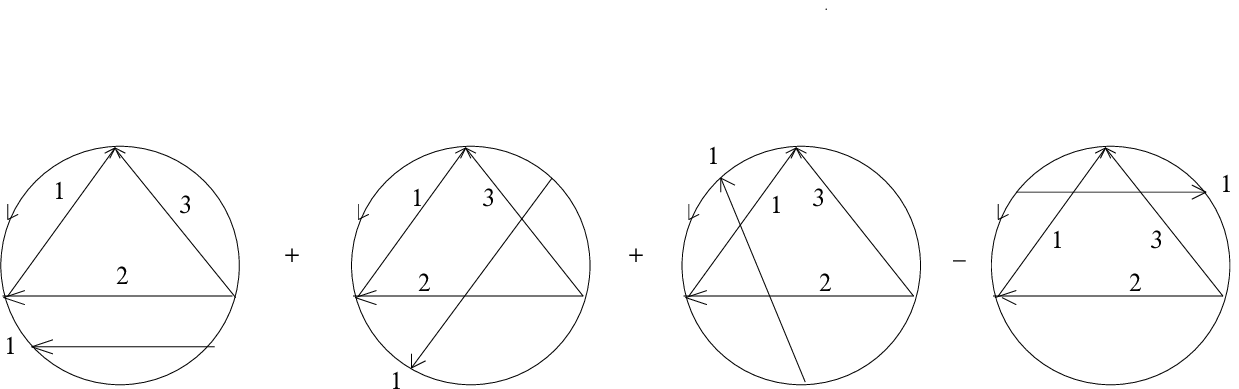}
\caption{}
\end{figure}
\begin{figure}[htbp]
\centering \psfig{file=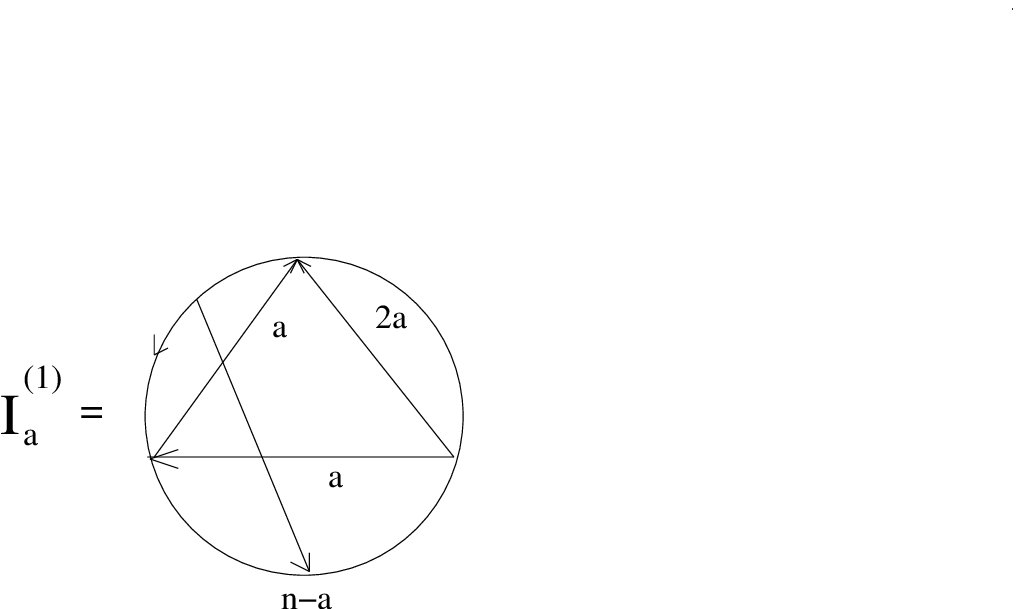}
\caption{}
\end{figure}
{\em Observation:\/} Let $a \in \{ 1, \dots ,[n/2]-1 \}$ be fixed
and let $I_a^{(1)}$ be the configuration shown in Figure 28.
This configuration is derived from the configuration $I$ in Figure 20
(see Remark 8).
Therefore, the configuration $I_a^{(1)}$ is a rigid braid-possible
configuration. But moreover, $I_a^{(1)}$ is a t-invariant
configuration and hence, $I_a^{(1)}$ {\em already defines a 1-cocycle
of degree 2\/}. In the next section, we will see how to generalize
$I_a^{(1)}$ in order to obtain non-trivial 1-cocycles of arbitrary degree
in a very simple way.
\subsection{One-cocycles of arbitrary degree}
For the degree $d>2$, other strata in $\Sigma^{(2)}$ also
impose conditions on $\Gamma$.

\begin{definition}
$\Gamma$ verifies the {\em tan-condition\/} if no
configuration $I$ in $\Gamma$ contains any sub-configuration
as shown in Figure 16.
\end{definition}

If $\Gamma$ verifies the {\em tan-condition\/}
then $\Gamma(S)$ is invariant under deformation of $S$ through points
in $\Sigma^{(1)}(tri) \cap \Sigma^{(1)}(tan)$.
The next definition is a straightforward generalization of Definition 14.

\begin{definition}
An {\em adjacent configuration\/} of a braid possible configuration $I$
is any configuration that is obtained in the following way:
one chooses an arrow among the $(d-1)$ that are not part of the triangle,
and one slides an endpoint of this arrow over exactly one endpoint of
another arrow (the latter may belong to the triangle).
All markings of arrows are preserved. The resulting configuration
has to be also braid possible.
\end{definition}

\begin{definition}
$\Gamma$ verifies the {\em tri-condition\/} if for each
(braid possible) configuration $\epsilon I$ in $\Gamma$, all adjacent
configurations of $I$ are also contained in $\Gamma$, with the same
coefficient $\epsilon$.
\end{definition}
If $\Gamma$ verifies the tri-condition, then $\Gamma(S)$
is invariant under deformation of $S$ through points in $\Sigma^{(1)}(tri)
\cap \Sigma^{(1)}(tri)$. The quadruple-condition for degree 2 is
generalized for arbitrary degree in the obvious way. Notice that if
$\Gamma$ satisfies the tri-condition then it satisfies automatically
the quadruple-condition.

Obviously, the strata of $\Sigma^{(2)}_3$ do not impose any condition
on $\Gamma$. The t-condition for degree 2 has to be
generalized in the following way.

\begin{definition}
The configurations on the same line in Figure 29 are called
{\em associated configurations\/} if besides the shown sub-configuration,
the rest of the configurations are identical. Notice that in the small
arcs of the circle there are no other endpoints of arrows.
\end{definition}
\begin{figure}[htbp]
\centering \psfig{file=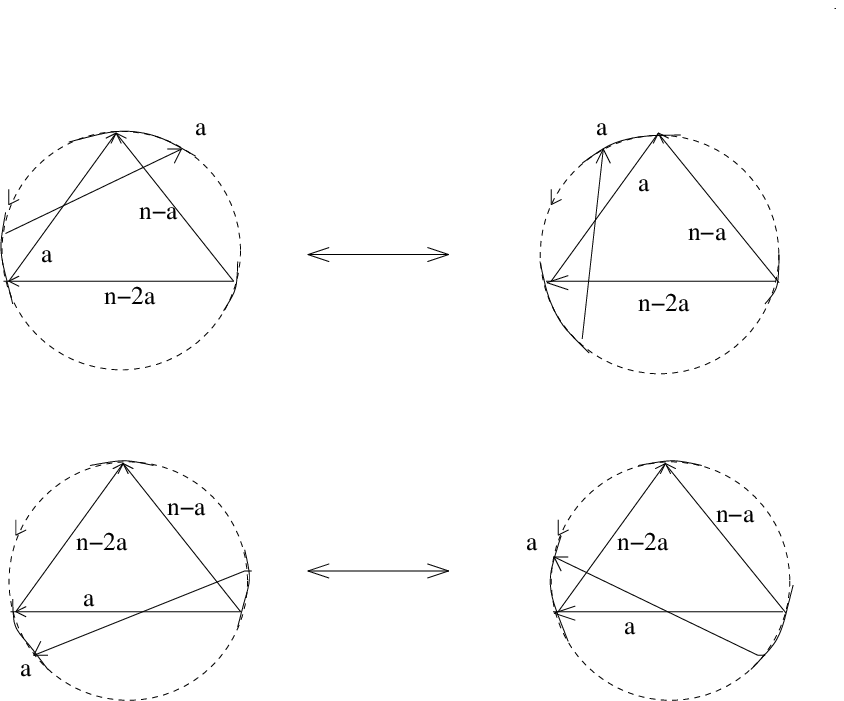}
\caption{}
\end{figure}
(Of course we extend all definitions to the configurations derived from
those of Figure 29 by the Remark 8).
If a configuration $I$ is different from all configurations in
Figure 29 (and of their derived configurations), then the associated
configuration will be the empty configuration.

\begin{definition}
$\Gamma$ satisfies the {\em t-condition\/} if each of its
configurations $\epsilon I$ occurs together with $-\epsilon I'$ where
$I'$ is the associated configuration.
\end{definition}

\begin{theorem}
 Let $\Gamma$ be a 1-cochain of degree $d$ that satisfies the
t, tri, tan-conditions. Then $\Gamma$ is a 1-cocycle of degree $d$.
\end{theorem}

{\em Proof:\/} It is completely similar to the proof of Theorem 4.$\Box$

\begin{remark}
Let us take a 1-cocycle of degree 1, and consider the sum of all
diagrams that are obtained from this 1-cocycle by adding one single arrow
with a new marking, in any possible position with respect to the triangle.
This sum is a 1-cocycle of degree 2. But this 1-cocycle is not interesting
because it is just the product of the 1-cocycle of degree 1 with
an invariant of degree 1 (see Proposition2).
In order to get a 1-cocycle $\Gamma$ which do not decompose into products
of 1-cocycles of lower degrees we need that {\em not all\/} possible
positions of an arrow with a fixed marking enter into $\Gamma$ with
the same coefficient.
\end{remark}
{\em Observation:\/} Let $\beta \in B_n$ and let $a \in \{ 1, \dots, n-1 \}$
be fixed. Then the sub-configurations shown in Figure 30 are {\em
locally rigid\/}, i.e. none of the two arrows can slide over the other
one in the small pictured arcs, because the resulting configuration
would not be braid possible. Moreover, they verify the
tan-condition.
Using this observation, it is easy to construct non-decomposing
1-cocycles of arbitrary degrees. We give an example in the following
proposition.
\begin{proposition}
 The configuration in Figure 31 defines a 1-cocycle of odd
degree $d$ if  $3a = n$.
\end{proposition}
{\em Proof:\/} Obviously, the whole configuration is rigid.
Moreover, no arrow with marking $n-a$ can become an arrow of the
triangle (by passing $\Sigma^{(2)}_2$), because it is always separated
from the triangle by an arrow with marking $a$. Evidently, no arrow with
marking $a$ can become an arrow of the triangle (by passing
$\Sigma^{(2)}_2$). No arrow can slide over the triangle, because $3a = n$. $\Box$
We generalize now the 1-cocycle $I^{(1)}_a$ from the end of the
previous section.

\begin{definition}
Let $d \in \mathbb{N}^*$ be odd and let $a \in \{ 1, 2, \dots, [n/2]-1 \}$
be fixed. The configuration $I_a^{(d)}$ is defined in Figure 32
(in the actual picture, we took $d=5$). The arrows with markings
$a$ and $n-a$ are alternating in the figure.
\end{definition}
\begin{proposition}
 $I^{(d)}_a$ defines a 1-cocycle of degree $d+1$.
\end{proposition}

{\em Proof:\/} It is completely analogous to the proof of Proposition 4.
$\Box$

\begin{definition}
Let $d \in \mathbb{N}^*$ be even, let $n$ be divisible by 3 and
let $a=n/3$. The 1-cochain $\Gamma_{n/3}^{(d)}$ is defined in
Figure 33 (in the actual picture, we took $d=4$).
\end{definition}
\begin{figure}[htbp]
\centering \psfig{file=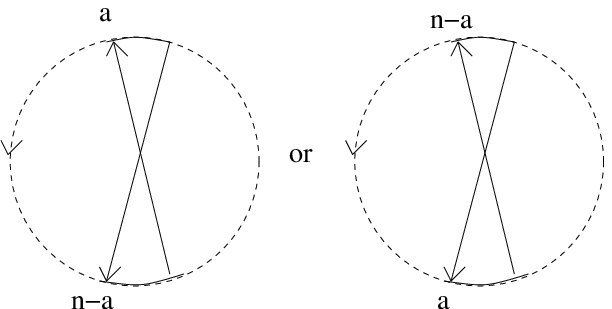}
\caption{}
\end{figure}
\begin{figure}[htbp]
\centering \psfig{file=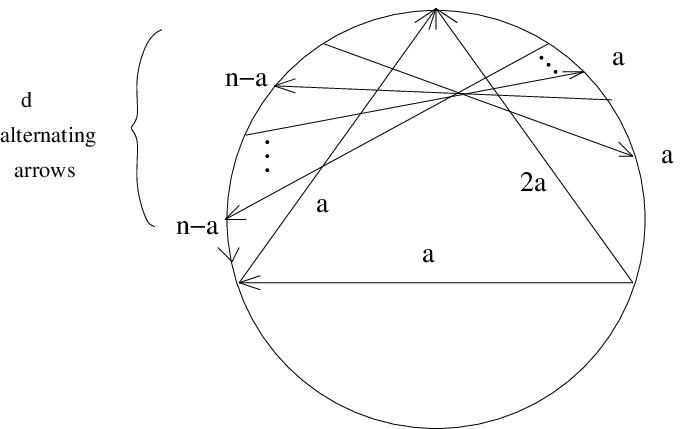}
\caption{}
\end{figure}
\begin{figure}[htbp]
\centering \psfig{file=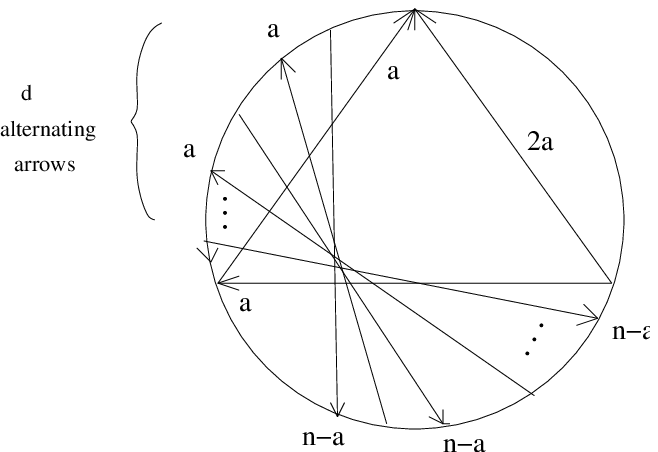}
\caption{}
\end{figure}
\begin{figure}[htbp]
\centering \psfig{file=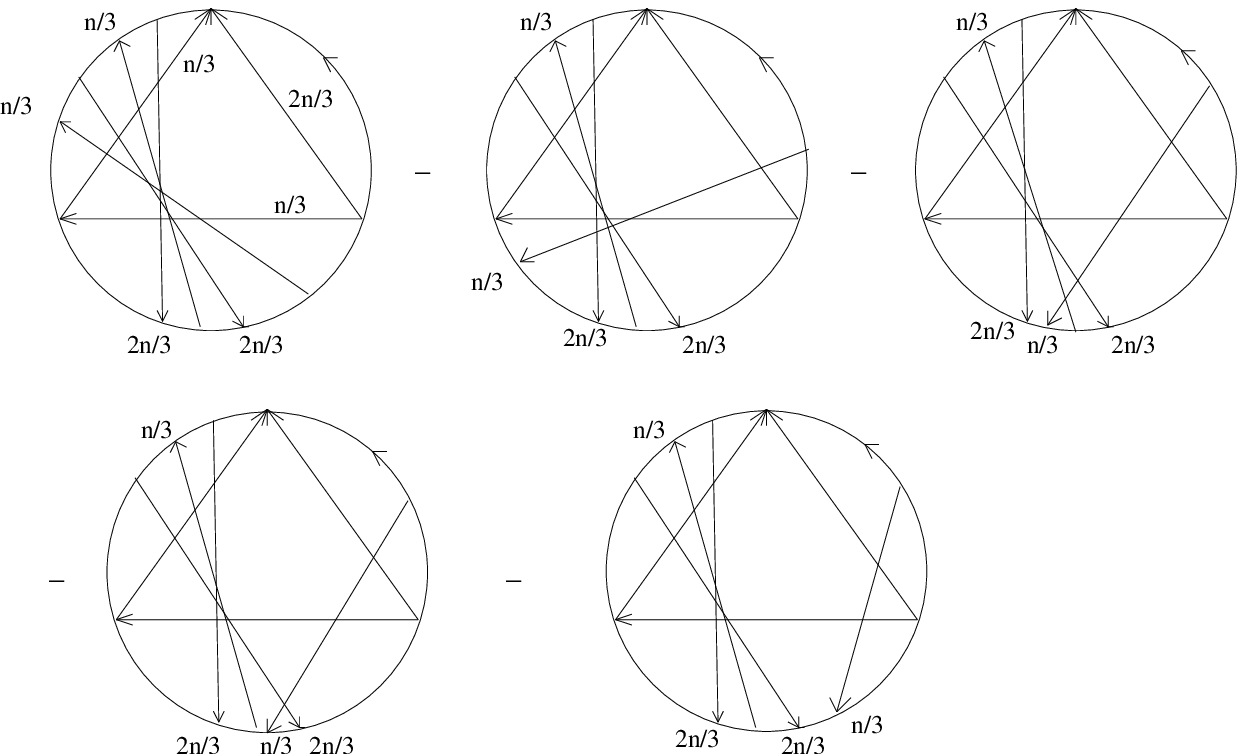}
\caption{}
\end{figure}
\begin{proposition}
 $\Gamma^{(d)}_{n/3}$ defines a 1-cocycle of degree $d+1$.
\end{proposition}

{\em Proof:\/} The first configuration in $\Gamma^{(d)}_{n/3}$ is
still rigid. But exactly one of the arrows with marking $n/3$
can be interchanged now with exactly one of the arrows in the triangle
(by passing $\Sigma^{(2)}_2$). The result is the second configuration
in the figure. This configuration is no longer rigid. The remaining
configurations in Figure 33 are just all the adjacent configurations.
We need that $3a=n$ in order to guarantee that in the last configuration
the arrow $n/3$ cannot slide further over some of the remaining
two vertices of the triangle.$\Box$

\begin{remark}
For example, the braid possible sub-configuration $II$ in Figure 20 with
$a = n/3$ is not contained in $\Gamma^{(d)}_{n/3}$. Such considerations
imply easily that in fact $\Gamma^{(d)}_{n/3}$ is not a product of
1-cocycles of lower degrees.

If $n$ is not divisible by 3 then we replace $\hat \beta$ by a
$3k$-cable, $k \in \mathbb{N}^*$, and take $a=kn$.

The mirror image $\Gamma^{(d)}_{n/3}!$ of $\Gamma^{(d)}_{n/3}$
(obtained by reversing all arrows, including those of the triangles,
and replacing all markings by their opposites) is of course also a
1-cocycle of degree~$d+1$.
\end{remark}

\begin{example}
Let $\beta = \sigma_1^{-1}\sigma_2^3 \in B_3$. Then,
$\Gamma_1^{(2)}!(rot(\hat \beta))=-1$. This shows that the
above 1-cocycles are not always trivial.
\end{example}
\begin{theorem}
 Let $K \hookrightarrow V$ be a closed braid which is a knot and let $\Gamma$
be a 1-dimensional cohomology class of degree $d$. Let $[rot(K)]$
be its canonical class. Then, $\Gamma([rot(K)])$ is a
$\mathbb{Z}$-valued finite type invariant of degree at most $d$ for
$K \hookrightarrow V$.
\end{theorem}

{\em Proof:\/} We have shown that $\Gamma([rot(K)])$ depends only
on the isotopy type of $K \hookrightarrow V$. The cocycle invariant
$\Gamma([rot(K)])$ is calculated as some sum $\sum_{s_i}$ over
triple points $s_i$ in $rot(K)$. Therefore, it suffices to
prove that this sum $\sum_{}$ for each triple point $s_i$ is of
finite type (even if it is not invariant).
If $\Gamma$ is of degree $d$ then $\sum_{s_i}$ depends only on the triple
point and of configurations of $d-1$ other crossings. This means that
in order to calculate a summand in $\sum_{s_i}$ we can switch all other
crossings besides the triple point and the fixed $d-1$ crossings. The
result will not change. This implies immediately that each $\sum_{s_i}$
is of degree $d$ (see \cite{PV} and also \cite{F01}).$\Box$

\begin{definition}
We call $\Gamma([rot(K)])$ a {\em 1-cocycle invariant of degree $d$.\/}
\end{definition}

\begin{example}
Let $K \hookrightarrow V$ be a closed 4-braid. Using Theorem 6,
Proposition 2, Proposition 3 and the examples of the next
section, one easily calculates that
\begin{displaymath}
\Gamma_{(2,1)^-}([rot(K)]) \equiv \Gamma_{(3,2)^+}([rot(K)])
\equiv W_1(K)-W_2(K)
\end{displaymath}
\end{example}
The 1-cocycle invariants have the following nice property, which can be
used to estimate from below the length of conjugacy classes of braids.
\vspace{0,5cm}

{\bf Proposition 1}
{\em Let the knot $K = \hat \beta \hookrightarrow V$ be a closed
$n$-braid and let $c(K)$ be its minimal crossing number, i.e. its
minimal word length in $B_n$. Then all 1-cocycle invariants
of degree $d$ vanish for\/}
\begin{displaymath}
d \geq c(K) + n^2 - n - 1.
\end{displaymath}

{\em Proof:\/}
Assume that the word length of $\beta$ is equal to $c(K)$.
We can represent $[rot(K)]$ by the following
isotopy which uses shorter braids.
\begin{displaymath}
\beta \to \Delta\Delta^{-1}\beta \to \Delta^{-1}\beta\Delta
\to \Delta^{-1}\Delta\beta' \to \beta' \to \Delta\Delta^{-1}\beta'
\to \Delta^{-1}\beta'\Delta \to \Delta^{-1}\Delta\beta \to \beta
\end{displaymath}

Here, $\beta'$ is the result of rotating $\beta$ by $\pi$ i.e.
each $\sigma_i^{\pm 1}$ is replaced by $\sigma_{n-i}^{\pm 1}$.
Obviously, $c(\Delta) = \frac{n(n-1)}{2}$. Thus, each Gauss diagram which
appears in the isotopy has no more than $c(K) + n^2 - n$ arrows. Indeed,
we create a couple of crossings by pushing $\Delta$ through $\beta$
only after having eliminated a couple of crossings before (see the
previous section). Therefore, for each diagram with a triple point
there are at most $c(K)+n^2-n-3$ other arows, and hence, each
summand in a 1-cocycle of degree $d$ is already zero if $d \geq
c(K) + n^2 - n - 1$.$\Box$
\subsection{The invariants are not functorial under cabling}

Bar-Natan, Thang Le, Dylan Thurston \cite{B-T-T} and independently S. Willerton \cite{Wi} have given a formula for the
Kontsevich integral of the cable of a knot in $\mathbb{R}^3$. This
formula has the following corollary:

\begin{corollary}(Bar-Natan, Thang Le, D. Thurston - Willerton) 
 Let $K, K'\hookrightarrow \mathbb{R}^3$ be knots which have the same
Vassiliev invariants up to a fixed degree $d$. Then, all (the same)
cables of $K$ and $K'$ have the same Vasiliev invariants up to degree
$d$.
\end{corollary}

In other words, it is useless to cable knots in purpose to
distinguish them by Vassiliev invariants.

Our 1-cocycle invariants of degree $d$ form a subset of all finite
type invariants of degree $d$ for knots in the solid torus.
For example, as already mentioned, there are no non-trivial 1-cocycle
invariants at all for closed 2-braids and there are no non-trivial
1-cocycle invariants of degree 1 for closed 3-braids. However, it turns out
that cabling is a usefull operation for the subset of 1-cocycle invariants
of degree $d$.
\begin{proposition}
The closed 2-braids $\hat \sigma_1$, $\hat \sigma_1^{-1}$ can
not be distinguished by any 1-cocycle invariants. However,
$Cab_2(\hat \sigma_1)$ and $Cab_2(\hat \sigma_1^{-1})$ can be
distinguished by a 1-cocycle invariant of degree 1.
\end{proposition}

{\em Proof:\/} $Cab_2(\hat \sigma_1)$ and $Cab_2(\hat \sigma_1^{-1})$
can be represented respectively by the 4-braids
$\beta = \sigma_3\sigma_2\sigma_1\sigma_2^2$ and
$\beta'= \sigma_3^{-1}\sigma_2\sigma_1^{-1}\sigma_2^{-2}$.
A calculation by hand shows:

\begin{displaymath}
\Gamma_{(2,1)^-}([rot(\hat \beta)]) =
\Gamma_{(3,2)^+}([rot(\hat \beta)]) =+1
\end{displaymath}

and

\begin{displaymath}
\Gamma_{(2,1)^-}([rot(\hat \beta')]) =
\Gamma_{(3,2)^+}([rot(\hat \beta')]) =-3
\end{displaymath}

$\Box$
\begin{remark}
Obviously, closed 2-braids are classified by the unique invariant
$W_1$ of degree 1. The braid $\hat \sigma_1$ is obtained from
$\hat \sigma_1^{-1}$ by multiplying $\sigma_1^{-1}$ with $\sigma_1^2$.
>From this, one easily concludes that for all $k \in \mathbb{Z}$,

\begin{displaymath}
\Gamma_{(2,1)^-}([rot(Cab_2(\sigma_1^{2k+1}))]) \equiv
\Gamma_{(3,2)^+}([rot(Cab_2(\sigma_1^{2k+1}))]) \equiv 1+4k
\end{displaymath}

Therefore, the unique non-trivial 1-cocycle of degree 1 for the
2-cable of 2-braids classifies also closed 2-braids.
\end{remark}

 Do all finite type invariants of closed braids and of local
knots arise as linear combinations of 1-cocycle invariants of
appropriate cables?

\subsection{Homological estimates for the number of braid relations in one-parameter families of closed braids}

Our 1-cocycles can be used in order to obtain information about
one-parameter families of closed braids which are knots.
Let $K$ be a closed n-braid which is a knot.
Let $S$ be a generic loop in $M(K)$.

\begin{definition}
The *-{\em length\/} $b([S])$ of $[S] \in H_1(M(K); \mathbb{Z})$ is
the minimal number of triple points in $S$ among all unions of generic
loops $S$ in $M(K)$ which represent $[S]$.
\end{definition}

\begin{theorem}
 Let $a, b \in \{ 0, 1, \dots, n \}$ with $a<b$. Then
\begin{displaymath}
b([S]) \geq 2 \sum_{(a,b)^+}{\vert \Gamma_{(a,b)^+}([S])\vert}
+2\sum_{(a,b)^-}{\vert \Gamma_{(a,b)^-}([S])\vert}
\end{displaymath}
\end{theorem}

{\em Proof:\/} Each $\Gamma_{(a,b)^+}$, $\Gamma_{(a,b)^-}$ is a
(in general non-trivial) 1-cocycle of degree 1 (see Proposition 3).
Each triple point in $S$ contributes by $\pm 1$ to the value of such a
cocycle. The inequality follows now from the relations (*) in
Proposition 3.$\Box$
\begin{example}
Let $\beta = \sigma_1\sigma_2^{-1}\sigma_3^{-1} \in B_4$.
For all $m \in \mathbb{Z}$, we have $b(m[rot(\hat \beta)])
\geq 4 \vert m \vert$.
This follows immediately from Example 2.
\end{example}
Theorem 7 does not contain any information in the case of closed
3-braids (because all 1-cocycles of degree 1 are trivial).
Therefore, we use the 1-cocycle $\Gamma$ of degree 2 from Example 2.
Let $\beta = \sigma_1\sigma_2^{-1} \in B_3$.
One has $\Gamma(rot(\hat \beta)) = +2$. It is easily seen that this
implies $\Gamma !(rot(\hat \beta)) = +2$ too, where $\Gamma !$
is the mirror image of $\Gamma$ (see Remark 10).
Thus, for all $m \in \mathbb{Z} \setminus 0$, $m(rot(\hat \beta))$
intersects both types of strata in $\Sigma^{(1)}(tri)$ (compare Figure 9). But the closure of the union of all strata of the same 
type in $\Sigma^{(1)}(tri)$  are
trivial cycles of codimension 1 in $M(\hat \beta)$ and hence,
$m(rot(\hat \beta))$ intersects each of these two strata in at
least two points.
It follows that $b(m[rot(\hat \beta)]) \geq 4$.
The canonical loop shows that this inequality is sharp for $m = \pm 1$.
\section{Character invariants}
In this section we introduce new easily calculable isotopy invariants for closed braids.

Remember, that a generic isotopy $\hat \beta_s , s \in [0 ,1]$, induces a tan-transverse homotopy $rot(\hat \beta_s) , s \in [0 ,1]$,
of the canonical loops.
The trace circles for different parameter s are in a natural one-to-one correspondence. Consequently, we can give {\em names \/}
$x_i$ to the circles of $TL\tilde(\hat \beta_0)$ and extend these names in a unique way on the whole family of trace circles.

Let $\{x_1 ,x_2 ,\dots \}$ be the set of named trace circles. Obviously, for each circle $x_i$ there is a well defined homological marking $h_i \in H_1(V)$.
Let $[x_i] \in H_1(T^2)$ be the homology class represented by the circle $x_i$ (with its natural orientation induced from the orientation of the trace graph).
\subsection{Character invariants of degree one}
In this section, we use the named cycles, i.e. the trace circles, in order to refine the
1-cocycles of degree 1 which were defined in Section 3.2.

Let K be a closed n-braid which is a knot.
Let $S \subset M(K)$ be a generic loop and let $X = \{ x_1, \dots, x_m \}$
be the corresponding set of cycles of named crossings. Let $x_{i_1}, x_{i_2}, x_{i_3} \in X$ be fixed.
We do not assume that they are necessarily different. Let $h_{i_1},
h_{i_2}, h_{i_3}$ be the corresponding homological markings.

\begin{definition}
A {\em character of degree 1 of $S$\/}, denoted by
\begin{displaymath}
C_{(h_{i_1}, h_{i_2})^\pm(x_{i_1}, x_{i_2}, x_{i_3})}(S)
\end{displaymath}
or sometimes shortly $C(S)$, is the algebraic intersection number of
$S$ with the strata $(h_{i_1}, h_{i_2})\pm$ in
$\Sigma^{(1)}(tri)$  and such that the crossings of the triple
point belong to the named cycles as shown in Figure 34.
We call the unordered set $\{ x_{i_1}, x_{i_2}, x_{i_3} \}$ the
{\em support\/} of the character $C$.
\end{definition}
\begin{figure}[htbp]
\centering \psfig{file=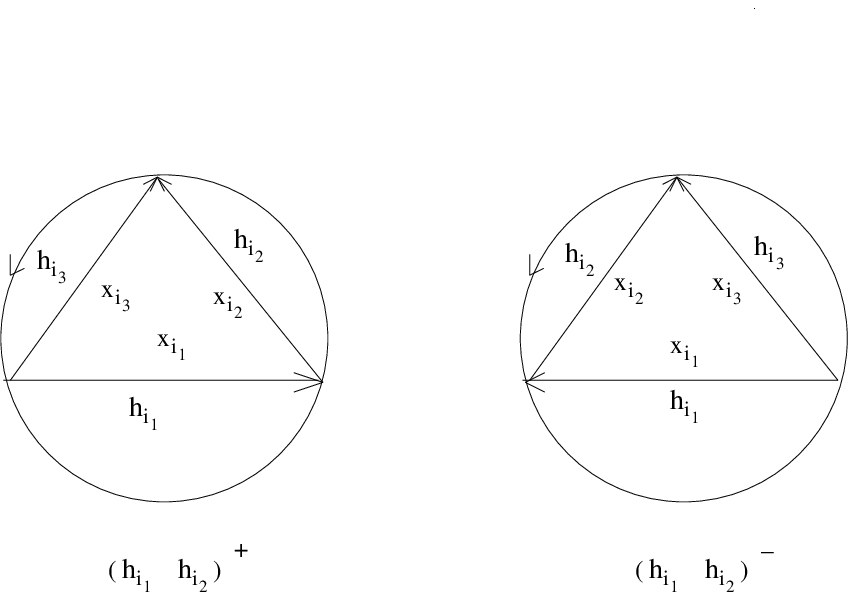}
\caption{}
\end{figure}
\begin{remark}
Evidently, in order to obtain a non-trivial intersection number we
need that $h_{i_3} = h_{i_1} + h_{i_2} - n$ in $(h_{i_1}, h_{i_2})+$
and that $h_{i_3} = h_{i_1} + h_{i_2}$ in $(h_{i_1}, h_{i_2})-$
(see Figure 10). It follows also that $x_{i_1} = x_{i_2} = x_{i_3}$
implies that necessarily $h_{i_1} = h_{i_2} = h_{i_3} =0$ or
$h_{i_1} = h_{i_2} = h_{i_3} = n$.
\end{remark}
Notice that for characters of degree 1 the relations (*) from
Proposition 3 are no longer valid. For example,
$C_{(h_{i_1}, h_{i_2})^+(x_{i_1}, x_{i_2}, x_{i_3})}(S)$
can be non-trivial for $h_{i_1} = h_{i_2}$ and even for
$x_{i_1} = x_{i_2}$.
\begin{theorem}
Let $S, S' \subset M(K)$ be generic loops and let
$\{ x_1, \dots, x_m \}$, $\{ x'_1, \dots, x'_{m'} \}$ be the corresponding
sets of named cycles. If $S$ and $S'$ are tan-transverse homotopic, then $m = m'$ and there is a bijection
$\sigma: \{ x_1, \dots, x_m \} \to \{ x'_1, \dots, x'_{m'} \}$ which
preserves the homological markings $h_i$ as well as the homology
classes $[x_i]$ and such that
\begin{displaymath}
C_{(h_{i}, h_{j})^\pm(x_{i_1}, x_{i_2}, x_{i_3})}(S) =
C_{(h_i, h_j)^\pm(\sigma (x_{i_1}), \sigma (x_{i_2}), \sigma (x_{i_3}))}(S')
\end{displaymath}
for all triples $(x_{i_1}, x_{i_2}, x_{i_3})$.
\end{theorem}

{\em Proof:\/} This is an immediate consequence of Lemma 5
 and the fact that the trace circles are isotopy invariants.$\Box$

In other words, characters $C(S)$ of degree one are invariants of
$[S]_{t-t}$.
\subsection{Character invariants of arbitrary degree for loops of closed braids}
We refine the results of the sections 3.3 and 3.4. in a
straightforward way.Let $S \subset M(K)$ be an oriented
generic loop and let $X = \{ x_1, \dots, x_m \}$ be the corresponding
set of named cycles.

Let $I$ be a configuration of degree $d$ (see Definition 11) and
let $(x_{i_1}, \dots, x_{i_{d+2}})$ be a fixed $(d+2)$-tuple of
elements in $X$ (not necessarily distinct). A {\em named
configuration\/} $I_{(x_{i_1}, \dots, x_{i_{d+2}}), \phi }$ is the
configuration $I$ together with a given bijection $\phi$ of
$(x_{i_1}, \dots, x_{i_{d+2}})$ with the $d+2$ arrows in $I$ and such
that $x_{i_1}, x_{i_2}, x_{i_3}$ are the arrows of the triangle
exactly as in the previous section. Of course, different bijections give
in general different named configurations.

A {\em character chain\/} $\Gamma$ of degree $d$ is a linear combination
$\Gamma = \sum_i{\epsilon_i I_i}, \epsilon_i \in \{ +1, -1, \}$
of named configurations $I_i$ of degree $d$ which are all defined
with the same $d+2$-tuple $(x_{i_1}, \dots, x_{i_{d+2}})$
where $x_{i_1}, x_{i_2}, x_{i_3}$ are the arrows of the triangle
as shown in Figure 34.
We refine the definitions and conditions of Sections 3.3. and
3.4. in the obvious way.
{\em Adjacent named configurations\/} are obtained by sliding the arrows
and preserving the names.

In {\em associated named configurations\/}, the two arrrows which
interchange have the same name. We show an example in Figure 35.
The {\em named quadruple and tri-conditions\/}
are now: $\Gamma$ contains all named adjacent configurations
and they enter $\Gamma$ with the same sign.

The {\em t-condition\/}
is now: $\Gamma$ contains all named associated configurations
and they enter $\Gamma$ with different signs.

The {\em tan-condition\/}
is now: no named configuration $I_i$ in $\Gamma$ contains any
sub-configuration as shown in Figure 36.
Notice that sub-configurations as shown in Figure 36 are allowed if the
arrows have different names $x_i \not= x_j$, even if they have the
same homological markings $h_i = h_j$.

\begin{theorem}
Let $S, S' \subset M(K)$ be generic loops and let
$\{ x_1, \dots, x_m \}$, $\{ x'_1, \dots, x'_m \}$ be the corresponding
sets of  named cycles. If $S$ and $S'$ are tan-transverse homotopic, then there is a bijection
\begin{displaymath}
\sigma: \{ x_1, \dots, x_m \} \to \{ x'_1, \dots, x'_m \}
\end{displaymath}
which preserves the homological markings $h_i$ as well as the
homology classes $[x_i]$ and such that
\begin{displaymath}
C_{(x_{i_1}, \dots, x_{i_{d+2}})}(S) =
C_{(\sigma (x_{i_1}), \dots, \sigma(x_{i_{d+2}}))}(S')
\end{displaymath}
for all character chains of degree $d$ which satisfy the named
quadruple, tri, t, tan-conditions.
\end{theorem}

{\em Proof:\/} This is completely analogous to the proofs of
Theorems 5 and 8.$\Box$

\begin{definition}
The character chains $C(S)$ of the above theorem are called
{\em characters of degree $d$\/}. They are invariants of
$[S]_{t-t}$.
\end{definition}\begin{example}
Here, all names $x_i$ are different.

The named configuration shown in Figure 37 defines a character of degree
5 for closed $n$-braids if $h_1 = h_2 = h_4 \leq [\frac {n-1}{2}]$.
Indeed, the configuration is braid possible and rigid (see Lemmas 6
and 7) and neither $x_2$ nor $x_4$ can become an arrow of the
triangle by passing $\Sigma^{(2)}_2$.

\end{example}
\begin{figure}[htbp]
\centering \psfig{file=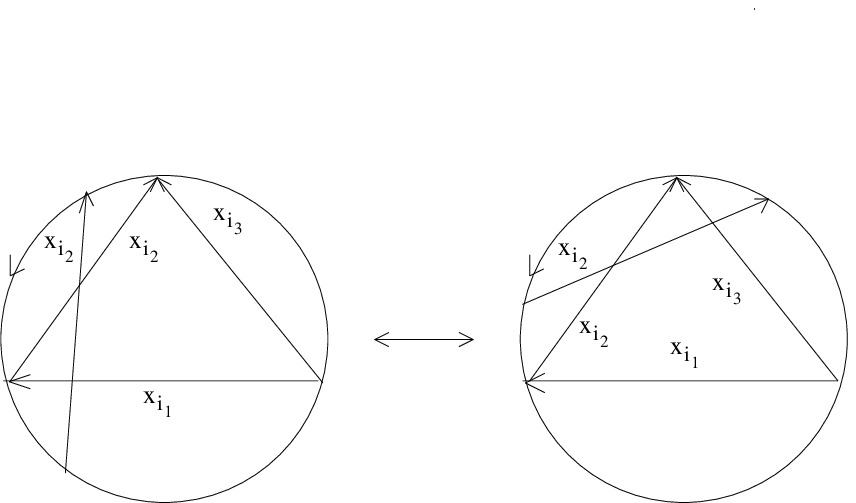}
\caption{}
\end{figure}
\begin{figure}[htbp]
\centering \psfig{file=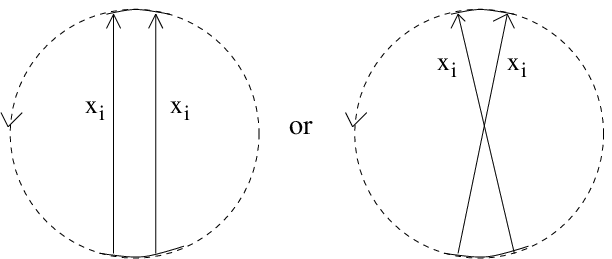}
\caption{}
\end{figure}
\begin{figure}[htbp]
\centering \psfig{file=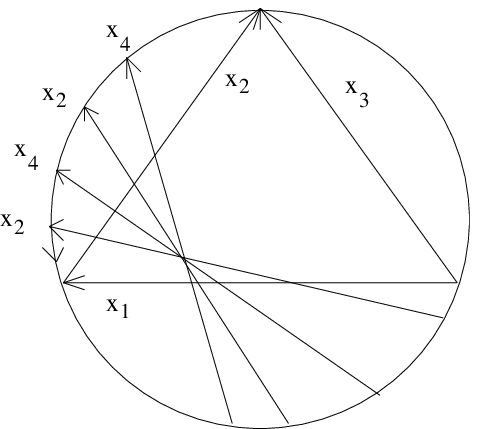}
\caption{}
\end{figure}
Let $h_i ,h_j$ and a type of triple point, e.g. $(h_i ,h_j)^+$, be fixed. It follows immediately from the definitions that

(**)            $\Gamma_{(h_i ,h_j)^+} = \sum_{(x_k ,x_l ,x_m)} C_{(h_i ,h_j)+}(x_k ,x_l ,x_m )$

where $h(x_k) = h_i$ and $h(x_l) = h_j$.

Hence, character invariants define splittings of one-cocycle invariants. However, the set of character invariants on the
right hand side of (**) is not an ordered set. Therefore we have to consider the set of character invariants as a local system on $M(\hat \beta)$.

It follows from Lemma 3  that the names $x_i$ are determined by their homological markings $h_i$, and that there are exactly n-1 trace circles.

However, this is in general no longer true in the case of multiples of the canonical loop.

Let $l \in \mathbb{N}$ be fixed and let $l rot(\hat \beta)$ be the loop which is defined by going $l$ times along the canonical loop.
Let $TL(l\hat\beta)$ denote the corresponding trace graph (the t-coordinate in the thickened torus covers now $l$ times the t-circle).

We will show in a simple example that the local system of character invariants is in general non trivial if $l > 1$.

\begin{example}
Let $\beta = \sigma_2\sigma_1^{-1}\sigma_2\sigma_1^{-1} \in B_3$.
We will write shortly $\beta = 2\bar12\bar1$.

The combinatorial canonical loop for $l = 2$ is given by the following sequence (where we write the names of the crossings just below the 
crossings).

$2\bar 12\bar 1 \hspace{0.3cm} \to \hspace{0.3cm}
\bar 1\bar 2\bar 12\bar 12(\bar 11)21 \hspace{0.3cm} \to \hspace{0.3cm}
\bar 1\bar 2\bar 12\bar 1221 \hspace{0.3cm} \to \hspace{0.3cm}
\bar 1\bar 2\bar 12\bar 121(\bar 121) \hspace{0.4cm} *_1\to$\\
\vspace{0.1cm}
\hspace*{0.5cm} $abcd \hspace{1.3cm}
zyxabcdxyz \hspace{1.4cm}
zyxabcyz \hspace{1.5cm}
zyxabcu_1u_1yz$\\

$\bar 1\bar 2\bar 12\bar 1(212)1\bar 2 \hspace{0.3cm} *_2\to
\bar 1\bar 2\bar 12(\bar 11)211\bar 2 \to
\bar 1\bar 2\bar 12211\bar 2 \to
\bar 1\bar 2\bar 121(\bar 121)1\bar 2 \hspace{0.3cm} *_3\to$
\vspace{0.1cm}
\hspace*{0.5cm} $zyxabcu_1zyu_1 \hspace{1.2cm}
                 zyxabzu_1cyu_1 \hspace{0.5cm}
                 zyxau_1cyu_1 \hspace{0.5cm}
                 zyxau_2u_2u_1cyu_1$\\

$\bar 1\bar 2\bar 1(212)1\bar 21\bar 2 \hspace{0.5cm} *_4\to
(\bar 1\bar 2\bar 1121)1\bar 21\bar 2 \hspace{0.3cm} \to \hspace{0.3cm}
1\bar 21\bar 2 \hspace{0.2cm} \to \hspace{0.2cm}
\bar 1\bar 2\bar 11\bar 21(\bar 212)1 \hspace{0.2cm} *_5\to$
\vspace{0.1cm}
\hspace*{0.2cm} $zyxau_2cu_1u_2yu_1 \hspace{0.9cm}
                 zyxcu_2au_1u_2yu_1 \hspace{0.6cm}
                 u_1u_2yu_1 \hspace{0.5cm}
                 z_1y_1x_1u_1u_2yu_1x_1y_1z_1$\\

$\bar 1\bar 2\bar 11\bar 2112(\bar 11) \hspace{0.3cm} \to \hspace{0.3cm}
\bar 1\bar 2\bar 11\bar 2112 \hspace{0.3cm} \to \hspace{0.3cm}
\bar 1\bar 2\bar 11\bar 21(121)\bar 1 \hspace{0.3cm}*_6\to \hspace{0.3cm}
\bar 1\bar 2\bar 11(\bar 212)12\bar 1$
\vspace{0.1cm}
$z_1y_1x_1u_1u_2yy_1x_1u_1z_1 \hspace{0.2cm}
z_1y_1x_1u_1u_2yy_1x_1 \hspace{0.2cm}
z_1y_1x_1u_1u_2yy_1x_1v_1v_1 \hspace{0.3cm}
z_1y_1x_1u_1u_2yv_1x_1y_1v_1$\\

$*_7\to
\bar 1\bar 2\bar 1112(\bar 11)2\bar 1 \to
\bar 1\bar 2\bar 11122\bar 1 \to
\bar 1\bar 2\bar 11(121)\bar 12\bar 1 \hspace{0.3cm} *_8\to
(\bar 1\bar 2\bar 1121)2\bar 12\bar 1$
\vspace{0.1cm}
$z_1y_1x_1u_1v_1yu_2x_1y_1v_1 \hspace{0.4cm}
                 z_1y_1x_1u_1v_1yy_1v_1 \hspace{0.3cm}
                 z_1y_1x_1u_1v_1yv_2v_2y_1v_1 \hspace{0.3cm}
                 z_1y_1x_1u_1v_2yv_1v_2y_1v_1$\\

$\to 2\bar 12\bar 1$\\
\hspace*{0.5cm} $v_1v_2y_1v_1$

\vspace{0.5cm}

We have the identifications: $d=x$, $b=z$, $x=c$, $y=u_2$, $z=a$,
$u_1=z_1$, $u_2=x_1$, $x_1=u_1$, $y_1=v_2$, $z_1=y$.
This gives us:

$2\bar 12\bar 1 \to 2\bar 12\bar 1$\\
\hspace*{0.5cm} $aacc  \hspace{0.6cm} v_1v_2v_2v_1$

together with the names $u_1=z_1=x_1=u_2=y$.

The second rotation gives us:

$2\bar 12\bar 1 \to 2\bar 12\bar 1 \to 2\bar 12\bar 1$\\
\hspace*{0.5cm} $aacc \hspace{0.4cm} v_1v_2v_2v_1 \hspace{0.4cm}
v_3v_4v_4v_3$

together with the identification $v_1=v_2$, $y_3=v_4$, $z_2=v_2$
and with the names $u_1=z_1=x_1=u_2=y$, $u_3=z_3=x_3=u_4=y_2$.
The monodromy (i.e. how the set of crossings is maped to itself after the rotation) implies now: $a=v_3=v_4=c$.

Therefore, we have exactly four named cycles: $a, v_1, u_1, u_3$ for

$2\bar 12\bar 1 \to 2\bar 12\bar 1$ with $l = 2$.

Consequently, we have 

$2\bar12\bar1 \to 2\bar12\bar1 \to 2\bar12\bar1$ with the names $aaaa \to v_1v_1v_1v_1 \to aaaa$.

 Hence, $rot(\hat \beta)$ 
acts by interchanging $a$ and $v_1$ (as well as $u_1$ and $u_3$). Consequently, the local system is non trivial in this example.
\end{example}
\subsection{The examples}
The following examples are calculated by Alexander Stoimenow using his program in c++. His program is available by request (see \cite{S}).

Let $\beta = \bar12\bar1^32^3 \in B_3$.
( $\hat \beta$ represents the knot $8_9$.)

We want to show that the link $\hat \beta  \cup$ (core of complementary solid torus) is not invertible in $S^3$.
This is equivalent to show that $\beta$ is not conjugate to $\beta_{inverse} = 2^3\bar1^32\bar1$ (compare \cite{F01}).

Because $\beta$ is a 3-braid , the homological markings are in $\{1 ,2\}$.

Character invariants of degree one for $l = 1$ and $l = 2$ do not distinguish $\hat \beta$ from $\hat \beta_{inverse}$.
However, for $l = 3$ we obtain three different named cycles  $x_1 ,x_2 ,x_3$ of homological marking 1 and three different named cycles 
$y_1 ,y_2 ,y_3$ of marking 2.

We consider the set of nine character invariants of degree one which are of the following form (see Figure 38) , where $i ,j \in \{1 ,2 ,3\}$.
\begin{figure}[htbp]
\centering \psfig{file=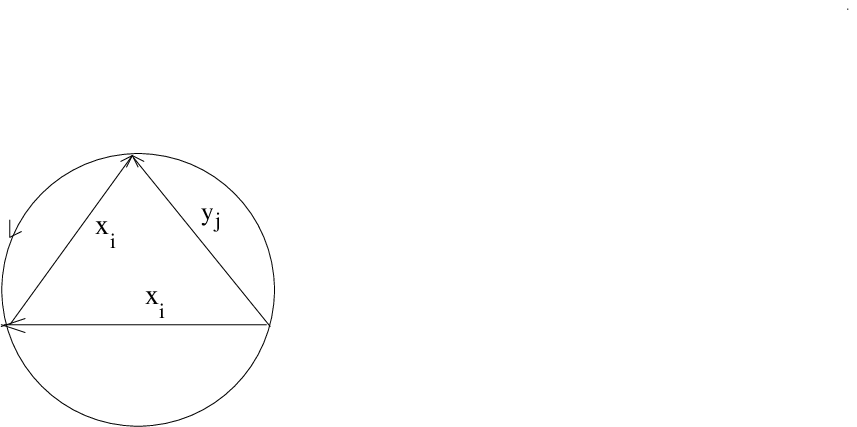}
\caption{}
\end{figure}
For $\hat \beta$  we obtain the set $\{-1 ,-1 ,-1 ,-1 ,-1 ,-1 ,2 ,2 ,2 \}$ and for $\hat \beta_{inverse}$ we obtain the set 
$\{1 ,1 ,1 ,1 ,1 ,1 ,-2 ,-2 ,-2\}$. Evidently, even without knowing the local system, there is no bijection of the trace circles for 
$\hat \beta$ and those for $\hat \beta_{inverse}$ which identifies the above sets.
Consequently, $\hat \beta$ and $\hat \beta_{inverse}$ are not isotopic.

The knot $9_5$ can be represented as a 8-braid with 33 crossings. Character invariants of degree one for $l = 2$ show that the braid is not 
invertible in the same way as in the previous example.

The knot $8_6$ can be represented as a 5-braid with 14 crossings. Character invariants of degree one for $l = 2 ,4 ,6$
show that it is not invertible as a 5-braid. (It does not work for for $l = 1 ,3 ,5$.)

The knot $8_{17}$ is not invertible as a 3-braid, which is shown with $l = 4$. (It does not work with $l = 1 ,2 ,3$.)

Let $b \in P_5$ be Bigelow's braid ( see \cite{Bi} ). It has trivial Burau representation. Let $s = \sigma_1\sigma_2\sigma_3\sigma_4$.
The braids s and bs have the same Burau representation. This is still true for their 2-cables, i.e. we replace each strand by two parallel 
strands.
Character invariants of degree one for $l = 2$ show that the (once positively half-twisted) 2-cables of the above braids are not conjugate, 
and, consequently, the braids s and bs are not conjugate either.
\subsection{A refinement of character invariants of degree one}
The number of triple points in a trace graph can only change by trihedron moves, as follows from Theorem 3 .

\begin{definition}
A {\em generalized trihedron \/} is a trihedron which might have other triple points on the edges.
\end{definition}

Figure 39 shows a tetrahedron move which transformes a trihedron into a generalized trihedron. The generalized trihedron has still exactly two vertices.
\begin{figure}[htbp]
\centering \psfig{file=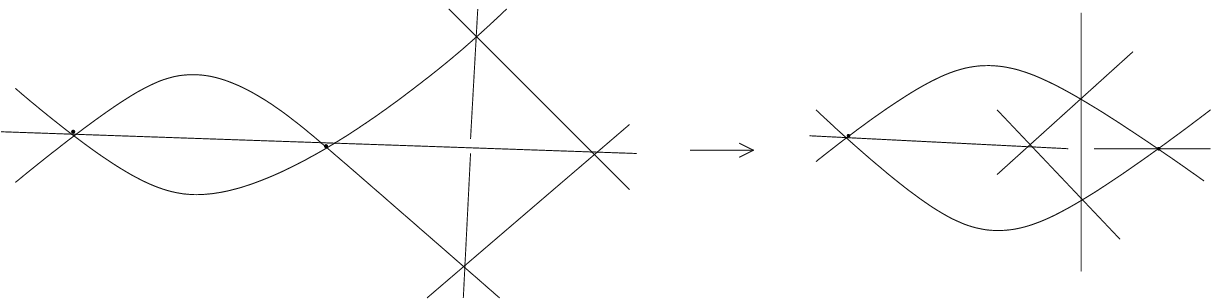}
\caption{}
\end{figure}
Evidently, the number of generalized trihedrons does not change under tetrahedron moves. Let $E$ be the set of all triple
points in the trace graph $TL(\hat \beta)$ which are not vertices of generalized trihedrons. The following lemma is an immediate consequence of Theorem 3.

\begin{lemma}
The set $E$, and, hence $card(E)$, is an isotopy invariant of closed braids $\hat \beta$.
\end{lemma}

Moreover, for each element of $E$ we have the additional structure defined before: type , sign, markings, names.

It follows from Theorem 3 and the geometric interpretation of generalized trihedrons in \cite{F-K06} that the two vertices of a generalized trihedron 
have always different signs. Consequently, character invariants of degree one count just the algebraic number of elements in $E$ which have a given type
and given names. But already the geometric number of such elements in $E$ is an invariant as shows Lemma 9.
\begin{definition}
Let $C^{+(-)}_{(h_i ,h_j)^{+(-)}}(x_k ,x_l ,x_m)$ be the number of all positive (respectively negative) triple points in $E$
of given type $(h_i ,h_j)^{+(-)}$ and with given names $x_k ,x_l ,x_m $. We call these the {\em positive (respectively negative) character invariants\/}.
\end{definition}  

The following proposition is now an immediate consequence of Lemma 9 and Definition 30.

\begin{proposition}
The positive and the negative character invariants are isotopy invariants of closed braids.
\end{proposition}

\begin{example}
Let us consider $\beta = \sigma_2\sigma_1^{-1} \in B_3$. Its trace graph $TL(\hat \beta)$ is shown in Figure 40.
One easily sees that it does not contain any generalized trihedrons. Consequently, all four triple points are in $E$.
There are exactly two names $x_1$ and $x_2$. They correspond to the homological markings
$h_1 = 1$ and $h_2 = 2$.

One easily calculates that two of the triple points are of type $(1 ,1)^-$ and they have different signs. The other two are of type 
$(2 ,2)+$ and they have different signs too. Consequently, all character invariants of degree one are zero.

However, we have $C^+_{(1 ,1)^-}(x_1 ,x_1 ,x_2) = C^-_{(1 ,1)^-}(x_1 ,x_1 ,x_2) = 1$, and\\
$C^+_{(2 ,2)^+}(x_2 ,x_2 ,x_1) = C^-_{(2,2)^+}(x_2 ,x_2 ,x_1) = 1$.Consequently, the positive and negative character invariants contain in this example more information than the character invariants of degree one (for $l = 1$).
\end{example}
There is not yet a computer program available in order to calculate these invariants in  more sophisticated examples.
\begin{figure}[htbp]
\centering \psfig{file=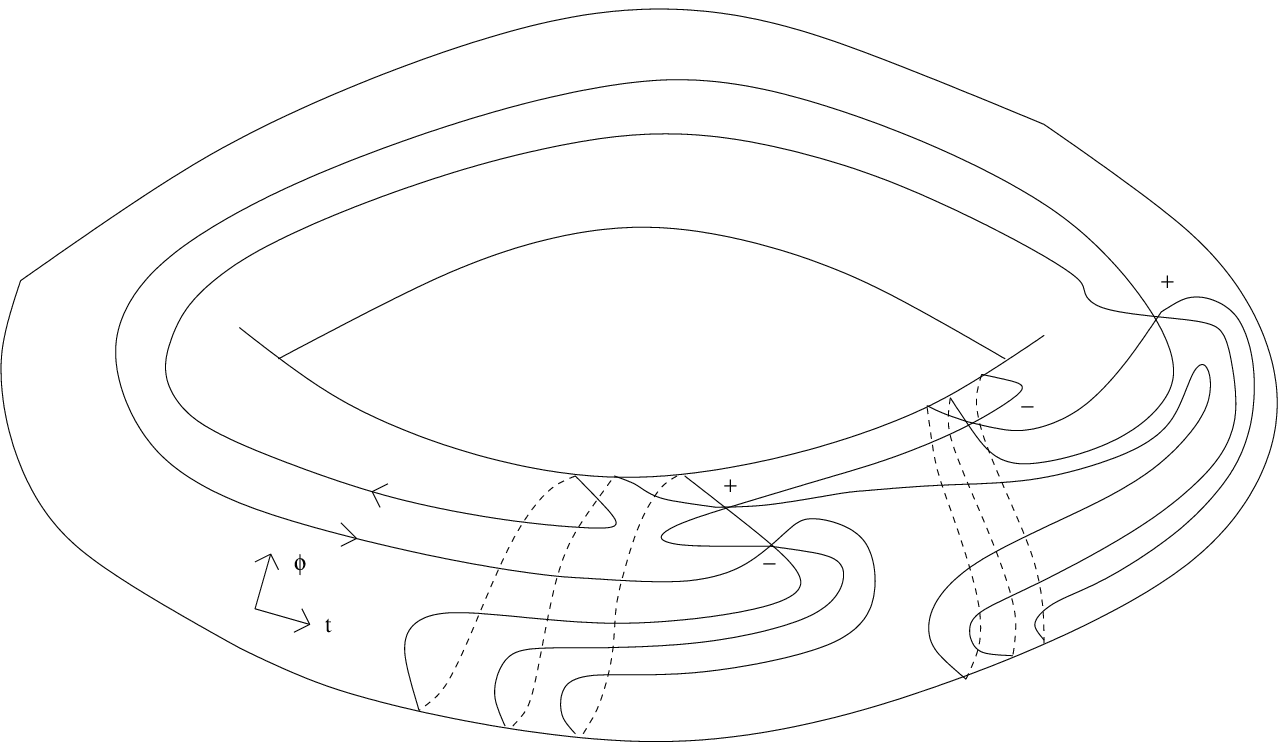}
\caption{}
\end{figure}
\subsection{Homotopical estimates for the number of braid relations in one parameter families of closed braids}

In section 3.6. we have used the 1-cocycles in order to estimate the
$*$-length $b([S])$ for classes $[S] \in H_1(M(K); \mathbb{Z})$.

Let $S \subset M(K)$ be a generic loop.

\begin{definition}
The $*$-{\em length\/} $b([S]_{t-t})$ of the tan-transvers homotopy class is the minimal number of triple points in $S$ amongst all
generic loops in $M(K)$ which represent $[S]_{t-t}$.
\end{definition}

\begin{theorem}
Let $\{ x_1, \dots, x_m \}$  be the set of essential named
cycles of $S$. Let $C_{(x_{i_1}, x_{i_2}, x_{i_3})}(S)$,
$\{ i_1, i_2, i_3 \} \subset \{ 1, \dots, m\}$ be the set of all
characters of degree one for $S$. Then

\begin{displaymath}
b([S]_{t-t}) \geq \sum_{\{ i_1,i_2,i_3 \} \subset \{ 1, \dots, m \}}
{\vert C_{(x_{i_1}, x_{i_2}, x_{i_3})}(S) \vert}
\end{displaymath}
\end{theorem}

{\em Proof:\/} Different triples $(x_{i_1}, x_{i_2}, x_{i_3})$
correspond to different strata of $\Sigma^{(1)}(tri)$.
Locally, each intersection index of $S$ with such a stratum is equal to
$\pm 1$. The result follows.$\Box$
\begin{example}
One easily calculates that for
$S=2rot(\hat {\sigma_2\sigma_1^{-1}})$ there are exactly eight
non-trivial characters of degree one. Each of them is equal to $\pm 1$.
This can be generalized straightforwardly for arbitrary $l$ with
$\vert l \vert \geq 2$. We obtain the following proposition:
\begin{proposition}
For all integers $l$ such that $\vert l \vert \geq 2$, we have:
\begin{displaymath}
b([l.rot(\hat {\sigma_2\sigma_1^{-1}})]_{t-t}) = 4\vert l \vert
\end{displaymath}
\end{proposition}

{\em Proof:\/} This follows immediately from Theorem 10 together
with a direct calculation which shows that

\begin{displaymath}
b([l.rot(\hat {\sigma_2\sigma_1^{-1}})]_{t-t}) \leq 4\vert l \vert
\end{displaymath}

$\Box$
\end{example}

Characters of higher degree can be used, of course, in the same way as
1-cocycles of higher degree were used to estimate $b([S])$ (see
section 3.6.).
\section{Character invariants of degree one for almost closed braids}

Let $K \hookrightarrow V = \mathbb{R}^3 \setminus z-axes$ be an oriented knot such that the restriction of $\phi$ to $K$
has exactly two critical points. In this case $K$ is called an {\em almost closed braid\/}. Necessarily, one of the critical points is a local maximum and 
the other is a local minimum. In analogy to the case of closed braids we consider almost closed braids up to isotopy through 
almost closed braids. 

We consider the (geometric) canonical loop $rot(K)$ and the corresponding trace graph $TL(K)$. $TL(K)$ is an oriented link with triple points and exactly four boundary points
 (corresponding to the oriented tangencies at the critical points of $\phi$).

There are more types of moves for  $TL(K)$ as in the case of closed braids (see \cite{F-K06}). But it turns out that only two of these additional moves are relevant 
for the construction of the invariants.

\begin{definition}
A {\em band move\/} is shown in Figure 41.
\end{definition}

\begin{definition}
A {\em unknot move \/} is shown in Figure 42.
\end{definition}
\begin{figure}[htbp]
\centering \psfig{file=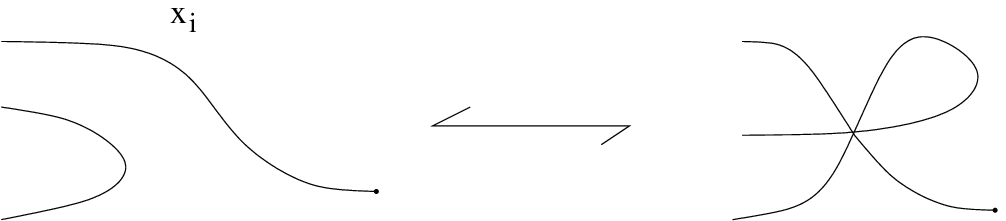}
\caption{}
\end{figure}
\begin{figure}[htbp]
\centering \psfig{file=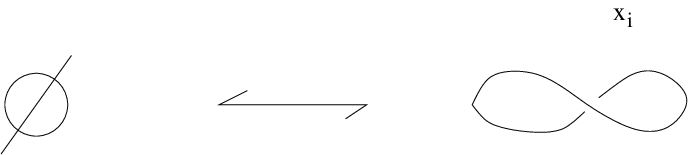}
\caption{}
\end{figure}
A band move corresponds to a branch which passes transversally through an ordinary cusp  and an unknot move corresponds to the case 
that the maximum and the minimum of $\phi$ have the same critical value. Notice, that there are no {\em extrem pair moves \/},
i.e. two maxima or two minima of $\phi$ with the same critical value  (for all this compare \cite{F-K06}).An extreme pair move induces a Morse modification of index 1 of the trace graph. Such modifications are difficult to control. We restrict ourselfs to almost closed braids in order to avoid them.

The unknot component $x_i$ of $TL(K)$ which was created by an unknot move, has always $[x_i] = 0$ in $H_1(T^2)$.
Notice, that in a band move there is always involved one component $x_i$ which has boundary.
Let $X = \{x_1 ,x_2 ,\dots \}$ be the set of all closed trace circles of $TL\tilde(K)$ which represent non-trivial homology classes in $H_1(T^2)$.

The above observations together with Theorem 8 and Theorem 1.10 from \cite{F-K06}  imply immediately the following theorem.

\begin{theorem}
Each character invariant of degree one with names $x_l ,x_k ,x_m \in X$ is an isotopy invariant for almost closed braids.
\end{theorem}

Evidently, we can apply the above theorem to $lrot(K)$ with arbitrary $l \in \mathbb{N}$ exactly as in the case of closed braids.

Laboratoire Emile Picard

Universit\'e Paul Sabatier

118 ,route de Narbonne 

31062 Toulouse Cedex 09, France

fiedler@picard.ups-tlse.fr
\end{document}